\documentclass[11pt,a4paper]{article}

\usepackage[T1]{fontenc}
\usepackage[utf8]{inputenc}
\usepackage{lmodern}
\usepackage[margin=25mm]{geometry}
\usepackage{microtype}
\usepackage{mathtools,amssymb,amsthm}
\usepackage{graphicx}
\usepackage{array,booktabs,multirow,threeparttable}
\usepackage{xcolor}
\usepackage{pifont}
\usepackage{algorithm,algpseudocode}
\usepackage{tikz}
\usetikzlibrary{arrows.meta,positioning,shapes.geometric}
\usepackage[font=small,labelfont=bf]{caption}
\usepackage{subcaption}
\usepackage{placeins,needspace}
\usepackage{authblk}
\usepackage[numbers,sort&compress]{natbib}
\usepackage[hyphens]{url}
\usepackage[hidelinks]{hyperref}
\hypersetup{
  pdftitle={Using Less for More: When Warm-Starting Accelerates Branch-and-Cut for Stochastic Programs},
  pdfauthor={Berend Markhorst, Joost Berkhout, Alessandro Zocca, Rob van der Mei}
}

\allowdisplaybreaks
\newtheorem{proposition}{Proposition}
\newcommand{\mymethod}{TULIP}
\newcommand{\cmark}{\ding{51}}

\newcommand{\sstpp}{SSTP\textsuperscript{+}}
\newcommand{\slrp}{SLR\textsuperscript{+}}
\newcommand{\fullmethodname}{``\textit{\textbf{T}hree-step warm-start framework \textbf{U}sed for solving \textbf{L}arge-scale stochastic mixed-\textbf{I}nteger \textbf{P}roblems}''}

\title{Using Less for More: When Warm-Starting Accelerates Branch-and-Cut for Stochastic Programs}
\author[1,2]{Berend Markhorst\thanks{Corresponding author: \href{mailto:berend.markhorst@cwi.nl}{berend.markhorst@cwi.nl}.}}
\author[2]{Joost Berkhout}
\author[2]{Alessandro Zocca}
\author[1,2]{Rob van der Mei}
\affil[1]{Stochastics Department, CWI, Science Park 123, 1098 XG Amsterdam, The Netherlands}
\affil[2]{Department of Mathematics, Vrije Universiteit Amsterdam, De Boelelaan 1105, 1081 HV Amsterdam, The Netherlands}

\date{}

\begin{document}
\maketitle

\begin{abstract}
    \setlength{\emergencystretch}{2em}
    Two-stage stochastic programs quickly become intractable as the number of scenarios grows. Motivated by this, we propose \mymethod, a modular and easy-to-implement three-step warm-start framework for two-stage stochastic (mixed-)integer programs with an exponential number of cuts separated during branch-and-cut. \mymethod\ (a)~builds a cheap \emph{surrogate} of the full problem by reducing the scenario set or by decoupling the two stages, (b)~solves it up to a first incumbent to collect the tight cuts separated along the way, and (c)~injects them to warm-start the original problem. In short: we use \textit{less} (a cheaper surrogate) for \textit{more} (the original problem). Using this modular setup, we propose four methods within this framework, each with a slightly different setting. Across four case studies, we show that this acceleration is governed by a single mechanism, the root cut loop, and we specify it through a closed-form equation. This \mymethod\ speedup model predicts a speedup when the time saved in the root cut loop exceeds the surrogate overhead. In our experiments, a \mymethod\ variant achieves mean speedups of up to 2.56, with gains increasing with the scenario count. In the remaining case studies, \mymethod\ provides little or no runtime benefit, which the \mymethod\ speedup model mostly explains through insufficient root cut loop savings compared to the surrogate overhead.

\par\medskip\noindent\textbf{Keywords:} Stochastic programming, Network design, Decomposition, Scenario reduction, Branch-and-cut, Warm start, Cut generation
\end{abstract}

\section{Introduction} \label{sec:introduction}
Stochastic programming (SP) is widespread across diverse fields in mathematical optimization, such as network design~\citep{ljubic_stochastic_2017}, supply chain planning~\citep{zahiri_multi-stage_2018}, energy markets~\citep{wallace_stochastic_2003}, airport operations~\citep{midler_stochastic_1969}, and inventory management~\citep{vaagen_modelling_2011}; see~\cite{birge_introduction_2011, klein_haneveld_stochastic_2020} for an overview. In SP problems, uncertain parameter outcomes correspond to branches in a scenario tree~\citep[Section 1.2.2]{king_modeling_2024}. A challenge in this field is that accurately approximating the underlying distributions often requires many scenarios, causing rapidly growing computational costs and making many problems intractable. This also happens with two-stage stochastic network design problems, whose formulations typically have an exponential number of valid inequalities, either in the first or second stage, that are added dynamically through branch-and-cut (B\&C)~\citep{ljubic_stochastic_2017}.

Several approaches make large stochastic programs easier to solve. Scenario reduction reduces the number of scenarios~\citep{chou_problem-driven_2023}, while decomposition uses the problem structure by splitting the stochastic program into smaller subproblems~\citep{song_adaptive_2015,pay_partition-based_2020}. These approaches are typically used to approximate or solve the original problem itself. A different idea is warm-starting, where information obtained from a related or smaller problem is reused to speed up the full solve~\citep{colombo_warm-start_2011}. Existing warm-start approaches for stochastic programs mainly transfer primal-dual or decomposition-specific information rather than problem-specific cuts. To our knowledge, no method uses such cuts collected from a cheaper surrogate to warm-start a generic branch-and-cut solve of the full stochastic problem.

To improve the tractability of exact solution methods for these problems, we propose a modular method with the idea of using less for more: \fullmethodname, or \mymethod\ in short. \mymethod\ consists of three steps. First, it constructs a cheaper surrogate of the full problem either by reducing the scenario set or by decoupling its two-stage structure. Second, the surrogate is solved until its first incumbent to generate cuts. Finally, these cuts, together with a (partial) primal start, warm-start the original problem with the full scenario tree.

In this work, we contribute to the existing literature as follows:
\begin{enumerate}
    \item We propose \mymethod, a modular and easy-to-implement three-step warm-start method for two-stage stochastic (mixed-)integer programs with an exponential number of valid inequalities in either stage.
    \item We derive a closed-form model for \mymethod's speedup that identifies the conditions under which constraint warm-starting accelerates branch-and-cut for two-stage stochastic programs and provide an upper bound on the attainable speedup using only information from the baseline branch-and-cut run without a warm start.
\end{enumerate}

Together, these contributions connect the design of a constraint-based warm-start method with an explanation of when such a method is computationally beneficial compared to a branch-and-cut. We evaluate \mymethod\ and this speedup model across four two-stage stochastic programs that differ in problem structure and in whether the exponentially large constraint sets appear in the first or second stage. For each case study, we compare \mymethod\ against regular branch-and-cut and, where applicable, compare surrogates based on scenario reduction and stage decoupling. We then explain the observed performance using our model in Equation~\eqref{eq:speedup-law} by comparing the speedup from fewer root cut loop iterations against the extra time required to build and solve the surrogate. This setup allows us to assess not only whether \mymethod\ accelerates branch-and-cut, but also \textit{why} it does or does not and how this depends on the choice of the surrogate. Using this setup, we address the following research questions:

\begin{enumerate}
    \item When does surrogate-based constraint warm-starting in \mymethod\ accelerate branch-and-cut?
    \item How well does the trade-off between reduced root cut loop time and surrogate overhead explain TULIP’s performance?
    \item How does the choice of surrogate construction affect \mymethod's performance?
\end{enumerate}

The remainder of this paper is structured as follows. Section~\ref{sec:literature} reviews the related literature. Section~\ref{sec:problem_class} introduces the problem setting. Section~\ref{sec:methodology} presents the TULIP framework and the \mymethod\ speedup model. Section~\ref{sec:computational_study} describes the computational study. Section~\ref{sec:results} reports and discusses the results. Finally, Section~\ref{sec:conclusion} summarizes our findings and outlines directions for future research.

\section{Related literature} \label{sec:literature}
Approaches for solving large-scale stochastic programs can be divided into heuristic and exact methods. Heuristic approaches reduce computational effort by simplifying or approximating the original problem, typically at the expense of an optimality guarantee. Exact approaches, in contrast, aim to improve computational tractability while retaining the guarantee of solving the original problem to optimality. \mymethod\ belongs to the latter class: although it uses a simplified surrogate to generate warm-start information, the original stochastic program with the full scenario set is then solved exactly by branch-and-cut.

The literature related to \mymethod\ can be partitioned into two groups. The first consists of methods for constructing computationally cheaper representations of stochastic programs, in particular through scenario reduction and decomposition. The second concerns warm-starting methods that reuse information from a related solve to accelerate the optimization process. We review these groups one by one and then position \mymethod\ relative to the closest approaches.

\subsection{Scenario reduction and decomposition} \label{sec:scen_red}
Large scenario sets can lead to intractable stochastic programs. One approach is therefore to replace the original scenario set with a smaller representation using scenario reduction methods, which can be either distribution-based~\citep{heitsch_scenario_2003, karuppiah_simple_2010} or problem-based~\citep{chou_problem-driven_2023, narum_problem-based_2024}. Distribution-based methods aim to preserve properties of the underlying probability distribution independently of the optimization problem, whereas problem-based methods use the problem structure to identify relevant scenarios, often yielding large reductions. The scenario reduction step in \mymethod\ is modular and, in principle, can use either type.

A related approach is Sample Average Approximation (SAA)~\citep{kleywegt_sample_2002}, which replaces the expected-value problem by one defined on a random sample of scenarios. Repeated samples can then be used to obtain statistical estimates of the solution quality. Both SAA and scenario reduction replace the full scenario set by a smaller one, thereby reducing the computational effort required to solve the stochastic program. Here, the resulting smaller problem serves as an approximation of the original problem. In \mymethod, however, the reduced problem is not used to provide the final solution. Instead, it functions as a cheaper surrogate from which information is transferred to an exact solve of the original problem.

Scenarios are also essential in decomposition methods for two-stage stochastic programs. Partition-based approaches~\citep{song_adaptive_2015, pay_partition-based_2020}, for example, group scenarios and initially represent each group through aggregate information, and refine the partition when more detail is needed. For the stochastic Steiner tree problem specifically, \cite{leitner_decomposition_2018} combine Benders decomposition and Lagrangian relaxation in a specialized exact algorithm. Such methods modify the solution procedure for the stochastic program itself. These methods change the way the stochastic program is solved. \mymethod\, in contrast, uses a simpler form of decoupling only to construct its surrogate: it considers the first-stage problem, the individual scenarios, or both independently. The resulting problems are used to generate warm-start information, while the branch-and-cut procedure for the original problem remains unchanged.

\subsection{Warm starts and cut reuse} \label{sec:warm_starts}
A second stream of literature accelerates optimization by transferring information obtained from another problem or an earlier solve. The transferred information can take different forms, including primal solutions, recourse information, and valid inequalities.

One possibility is to transfer a primal solution. \cite{colombo_warm-start_2011}, and subsequently \cite{colombo_decomposition-based_2013}, warm-start large stochastic programs solved by interior-point methods using a solution obtained from a smaller version of the original problem. Their results show that information from a reduced stochastic problem can reduce both runtime and the number of iterations. Although this idea is close to \mymethod\ in its use of a cheaper related problem, \mymethod\ transfers valid inequalities rather than a primal-dual warm-start point, and targets branch-and-cut rather than interior-point optimization.

More recently, \cite{kothari2026accelerating} study how information can be reused when solving a sequence of related SAA problems with Benders decomposition. They reuse dual solutions from earlier runs to generate Benders cuts for later ones, and they also study which cuts should be added to the master problem at the start. This is related to \mymethod\ because both approaches reuse information from one stochastic problem to speed up another. The main difference is that they work with a sequence of SAA replications and Benders cuts, whereas \mymethod\ transfers cuts from a cheaper surrogate to an exact branch-and-cut solve of the full problem.

Information from decomposition can also be reused to improve later optimization. For example, \cite{angulo_improving_2016} improve the integer L-shaped method~\citep{laporte_integer_1993} by switching between relaxations and mixed-integer subproblems when computing second-stage costs. This shows that information from cheaper intermediate problems can reduce the work needed by the exact algorithm. In \mymethod, the same idea is applied more directly: we collect cuts during the surrogate solve and add them to the original model before branch-and-cut starts.

Closer to our setting, valid inequalities generated in an auxiliary phase have been reused directly in branch-and-cut. For the deterministic Steiner tree problem, \cite{de2001dual} start branch-and-cut with active cuts obtained from a dual-ascent procedure~\cite{wong_dual_1984}. Similarly, \cite{fischetti2017thinning} generate connectivity cuts during a local-branching phase and use them to start the final branch-and-cut solve. These methods demonstrate that cuts generated cheaply before the main solve can be computationally valuable. In these approaches, an auxiliary method is first applied to the same deterministic problem, and the resulting cuts are then used to warm-start its final branch-and-cut solve. \mymethod\ instead generates cuts from a cheaper surrogate and transfers them to the original two-stage stochastic problem.

\subsection{Positioning of \mymethod} \label{sec:positioning_tulip}
The preceding literature covers three relevant ideas: (i) large stochastic programs can be made cheaper through scenario reduction or decomposition; (ii) information from a smaller problem can warm-start a larger solve; and (iii) cuts generated in an auxiliary phase can accelerate branch-and-cut. A gap remains at their intersection. Existing scenario reduction methods mainly use a smaller problem to approximate the original stochastic program, while existing warm-start and information-reuse methods transfer primal-dual or decomposition-specific information between related solves. Prior work demonstrates that information generated in auxiliary or reduced problems can be computationally useful, but does not identify when the cost of generating and transferring that information is outweighed by the resulting savings in the exact solve. In particular, it remains unclear when cuts collected from a surrogate stochastic problem can accelerate an exact branch-and-cut solve of the original problem.

\mymethod\ addresses this gap by solving a cheaper surrogate only long enough to collect useful cuts and then injecting these cuts into the original stochastic program before solving it exactly by branch-and-cut. The framework is modular in how the surrogate is constructed: in our experiments, it can be obtained via scenario reduction or stage decomposition. The contribution is therefore not a new scenario reduction or decomposition method, but a cut-transfer framework and an analysis of the computational conditions under which this transfer is useful. Table~\ref{tab:literature} summarizes this positioning.

\begin{table}[ht!]
\centering
\caption{Positioning of \mymethod\ against the most closely related approaches. \cmark\ indicates that the property holds.}
\label{tab:literature}
\footnotesize
\begin{threeparttable}
\begin{tabular}{@{}>{\raggedright\arraybackslash}p{4.1cm}*{4}{>{\centering\arraybackslash}p{1.6cm}}@{}}
\toprule
\textbf{Method} &
\textbf{Two-stage stochastic} &
\textbf{Cut transfer} &
\textbf{Generic branch-and-cut} &
\textbf{Payoff characterization} \\
\midrule
SAA~\cite{kleywegt_sample_2002} & \cmark & & & \\
Interior-point warm start~\cite{colombo_warm-start_2011} & \cmark & & & \\
Integer L-shaped~\cite{angulo_improving_2016} & \cmark & & & \\
Partition decomposition~\cite{song_adaptive_2015, pay_partition-based_2020} & \cmark & & & \\
Dual-ascent Steiner tree~\cite{de2001dual} & & \cmark & \cmark & \\
Local-branching Steiner tree~\cite{fischetti2017thinning} & & \cmark & \cmark & \\
Benders+Lagrangian Steiner tree~\cite{leitner_decomposition_2018} & \cmark & & & \\
Accelerated Benders~\cite{kothari2026accelerating} & \cmark & \cmark & & \\
\addlinespace
\textbf{\mymethod\ (ours)} & \cmark & \cmark & \cmark & \cmark \\
\bottomrule
\end{tabular}
\end{threeparttable}
\end{table}

\section{Problem setting} \label{sec:problem_class}
In this work, we focus on solving the following two-stage stochastic mixed-integer program to optimality:
\begin{equation} \label{eq:sp}
	\min_{x \in \mathcal{X}} \left\{ c^\top x + \mathbb{E}^{\mathbb{P}}\!\left[ v(x, \xi) \right] \right\},
\end{equation}
in which $x$ are the first-stage decisions leading to a direct cost of $c^\top x$. The feasible region for $x$ is
\begin{equation}\label{eq:first-stage}
    \mathcal{X} := \Bigl\{\, x \in \mathcal{Z}_1 \;\Big|\; A x \geq b,\;\;
    (a_j)^\top x \geq d_j \;\; \forall j \in \mathcal{C}^{(0)} \,\Bigr\},
\end{equation}
where $\mathcal{Z}_1$ is a mixed-integer domain. In the second stage, we incur an expected recourse cost of $\mathbb{E}^{\mathbb{P}}\!\left[ v(x, \xi) \right]$, where $\xi$ denotes the uncertain second-stage data that is only revealed after decision $x$ is taken. We represent the uncertainty by a finite set of scenarios $\mathcal{S}$: each scenario $s \in \mathcal{S}$ corresponds to a realization $\xi^{(s)}$ of $\xi$ that occurs with probability $p^{(s)} > 0$, with $\sum_{s \in \mathcal{S}} p^{(s)} = 1$. The distribution $\mathbb{P}$ is therefore supported on $\mathcal{S}$, so that the expected recourse cost corresponds to $\mathbb{E}^{\mathbb{P}}\!\left[ v(x, \xi) \right] = \sum_{s \in \mathcal{S}} p^{(s)}\, v(x, \xi^{(s)})$. For each scenario $s$, the feasible region for recourse decisions $y$ is
\begin{equation}\label{eq:second-stage-feasible-region}
	\mathcal{Y}^{(s)}(x) := \Bigl\{\, y \in \mathcal{Z}_2 \;\Big|\;
	W^{(s)}\, y \geq h^{(s)} - T^{(s)}\, x,\;\;
	(w^{(s)}_j)^\top y \geq e^{(s)}_j \;\; \forall j \in \mathcal{C}^{(s)} \,\Bigr\},
\end{equation}
where $\mathcal{Z}_2$ is a mixed-integer domain.
We assume relatively complete recourse, i.e., $\mathcal{Y}^{(s)}(x) \neq \emptyset$ for every $x \in \mathcal{X}$ and $s \in \mathcal{S}$. Based on scenario-based recourse optimization, the second-stage cost function is
\begin{equation}\label{eq:second-stage}
	v(x, \xi^{(s)}) := \min_{y \in \mathcal{Y}^{(s)}(x)} \Bigl\{\, (q^{(s)})^\top y \,\Bigr\}.
\end{equation}
Part of the constraints in $\mathcal{X}$ and $\mathcal{Y}^{(s)}(x)$ for every $s \in \mathcal{S}$ are collected in cut sets $\mathcal{C}^{(0)}$ and $\mathcal{C}^{(s)}$, respectively. In \eqref{eq:first-stage}--\eqref{eq:second-stage-feasible-region},
$\mathcal{C}^{(0)}$ and $\mathcal{C}^{(s)}$ denote index sets. With a slight abuse of notation, we also use them for the corresponding families of cuts. \mymethod\ focuses on problems in which the sizes of $\mathcal{C}^{(0)}$ and/or $\mathcal{C}^{(s)}$ for $s \in \mathcal{S}$ are exponential in the problem size. This is, for example, the case for connectivity constraints in network-design problems and subtour-elimination constraints in routing problems; see also Section~\ref{sec:computational_study}. Therefore, separating cuts from $\mathcal{C}^{(0)}$ and $\mathcal{C}^{(s)}$ for $s \in \mathcal{S}$ during branch-and-cut is required in practice. In contrast to these dynamically added cuts, the static constraints $Ax \ge b$ and $W^{(s)}y \ge h^{(s)} - T^{(s)}x$ from $\mathcal X$ and $\mathcal Y^{(s)}(x)$, respectively, are directly added to the initial formulation.

\section{The \mymethod\ framework} \label{sec:methodology}
This section presents \mymethod, a three-step warm-start framework for two-stage stochastic (mixed-)integer programs with an exponential number of cuts in either stage. \mymethod's possible variants are detailed in Section~\ref{sec:algorithm} together with a specification of its exactness. Section~\ref{sec:accounting_model} presents a model to predict the effectiveness of a \mymethod\ variant.

\subsection{\mymethod\ and its variants} \label{sec:algorithm}

\mymethod\ exploits the idea that cuts separated for a simpler surrogate representation are useful for the original problem while computationally cheaper to obtain. To make this concrete, it helps to view the two-stage stochastic program from~\eqref{eq:sp} as consisting of one first-stage block, indexed by $0$, and one recourse block for each scenario $s \in \mathcal{S}$. Then, a surrogate representation may be composed of one or more of these blocks, and the cuts separated for this surrogate can be reused in the original problem. We formalize this idea as follows.

Let $\mathcal B:=\{0\}\cup\mathcal S$ denote the set containing the first-stage block and all recourse blocks. For each $b \in \mathcal{B}$, $\mathcal{C}^{(b)}$ denotes the corresponding cut family. For $B \subseteq \mathcal{B}$, we then write $\mathcal{I} (B,\widehat{\mathcal{C}} )$ for the surrogate for \eqref{eq:sp} obtained by ``keeping'' the blocks in $B$ and imposing cuts $\widehat{\mathcal{C}} \subseteq \cup_{b \in B} \mathcal{C}^{(b)}$ upfront. Imposing the cuts $\widehat C$ upfront corresponds to warm-starting the surrogate with these cuts. When solving $\mathcal{I} (B,\widehat{\mathcal{C}} )$, the remaining cuts in $\cup_{b \in B} \mathcal{C}^{(b)}\setminus\widehat{\mathcal{C}}$ are separated dynamically during branch-and-cut. If the first-stage block is not considered, i.e., $0 \not\in B$, solving $\mathcal{I}(B,\widehat{\mathcal{C}})$ requires that the choice of $x$ be specified upfront. This choice is problem-specific, and Section~\ref{sec:computational_study} and the Appendix detail how to do it for the case studies we considered. For a reduced scenario set $\mathcal S'\subseteq\mathcal S$, let
$\tilde p^{(s)}>0$, $s\in\mathcal S'$, denote the probabilities assigned to the retained scenarios, with $\sum_{s\in\mathcal S'}\tilde p^{(s)}=1$.
The following surrogates $I(B,\widehat C)$ arise for different choices of $B\subseteq\mathcal B$, with $\widehat C\subseteq\bigcup_{b\in B}\mathcal C^{(b)}$ imposed upfront and $\mathcal S'\subseteq\mathcal S$:
\[
    \mathcal{I}(B,\widehat{\mathcal{C}})=
    \begin{cases}
        \displaystyle \min_{x\in\mathcal{X}} \left\{c^\top x\right\},
        & B=\{0\},\\[4pt]
        \displaystyle \sum_{s\in\mathcal S'} \tilde p^{(s)}
        v(x,\xi^{(s)}) \mbox{ for given } x,
        & B=\mathcal{S'},\\[4pt]
        \displaystyle \min_{x\in\mathcal X} \left\{ c^\top x + \sum_{s\in\mathcal S'} \tilde p^{(s)} v(x,\xi^{(s)}) \right\},
        & B=\{0\}\cup\mathcal{S'}.
    \end{cases}
\]
Note that the original problem as specified in \eqref{eq:sp} corresponds to $\mathcal{I} \bigl(\mathcal{B}, \cup_{b \in \mathcal{B}} \mathcal{C}^{(b)} \bigr)$, whereas $\mathcal{I}(\mathcal{B},\emptyset)$ denotes the same problem solved from scratch, with no cuts imposed upfront and all cuts separated dynamically during the solve. We refer to this latter standard approach as regular branch-and-cut.

Depending on the problem, different surrogates may be effective for finding cuts that accelerate the solution of the original problem. Moreover, different surrogates may emphasize different structural aspects and thus may identify different useful cuts for a single problem. Therefore, \mymethod\ allows the usage of multiple surrogates, each defined by a block set $B$. Let $\mathcal{G}$ denote a collection of block sets used. For illustration purposes in numerical experiments, this paper considers four TULIP variants corresponding to four different groupings of the problem blocks:
\[
\begin{aligned}
\text{TULIP-R(eduction)}: \qquad
    \mathcal{G}
    &=
    \left\{
        \{0\}\cup\mathcal{S}'
    \right\},\\
\text{TULIP-F(irst stage block)}: \qquad
    \mathcal{G}
    &=
    \left\{
        \{0\}
    \right\},\\
\text{TULIP-S(econd stage blocks)}: \qquad
    \mathcal{G}
    &=
    \left\{
        \{s\}:s\in\mathcal{S}
    \right\},\\
\text{TULIP-C(ombining all blocks)}: \qquad
    \mathcal{G}
    &=
    \left\{
        \{0\}
    \right\}
    \cup
    \left\{
        \{s\}:s\in\mathcal{S}
    \right\}.
\end{aligned}
\]

Algorithm~\ref{alg:tulip} details the TULIP steps using the block notation introduced above. In Step~1, a collection $\mathcal{G}$ of block sets is selected, where each $B\in\mathcal{G}$ defines one surrogate solve. Step~2 solves each surrogate until its first incumbent and collects the cuts that are tight at that incumbent. Finally, Step~3 solves the full problem with the collected cuts imposed upfront and, whenever compatible, incumbent information from the surrogate solves supplied as a primal start. Except where stated otherwise, a solve continues to optimality.

\begin{algorithm}[ht]
\caption{TULIP}
\label{alg:tulip}
\begin{algorithmic}[1]

\Statex \textbf{Step 1: Construct the surrogate problem(s).}

\State Select a collection $\mathcal{G}$ of block sets
       $B \subseteq \mathcal{B}$.

\Statex
\Statex \textbf{Step 2: Collect cuts from the surrogate problem(s).}

\State $\overline{\mathcal{C}} \gets \emptyset$

\For{$B \in \mathcal{G}$}
    \State Solve $\mathcal{I}(B,\emptyset)$ by branch-and-cut
           until the first incumbent.
    \State Let $\widehat{\mathcal{C}}_B \subseteq {\cup_{b \in B} \mathcal{C}^{(b)}}$
           denote the cuts separated during the solve.
    \State $\overline{\mathcal{C}}
           \gets
           \overline{\mathcal{C}}
           \cup
           \left\{
               c \in \widehat{\mathcal{C}}_B :
               c \text{ is tight at the returned incumbent}
           \right\}$.
\EndFor

\Statex
\Statex \textbf{Step 3: Solve the full problem.}

\State Solve
       $\mathcal{I}(\mathcal{B},\overline{\mathcal{C}})$
       by branch-and-cut until optimality, using compatible
       incumbent information from Step~2 as a primal start.

\end{algorithmic}
\end{algorithm}

Since $\overline{\mathcal{C}}$ in Algorithm~\ref{alg:tulip} is only composed of original cuts $\cup_{b \in \mathcal{B}} \mathcal{C}^{(b)}$ from \eqref{eq:sp} by construction, solving $\mathcal{I}(\mathcal{B},\overline{\mathcal{C}})$ in Step~3 is guaranteed to return an optimal solution to the original problem $\mathcal{I}(\mathcal{B},\cup_{b \in \mathcal{B}} \mathcal{C}^{(b)})$.

\begin{proposition}[Exactness of \mymethod]
\label{prop:exactness}
Algorithm~\ref{alg:tulip} returns an optimal solution to the full problem~\eqref{eq:sp}.
\end{proposition}
Importantly, exactness does not require the transferred cuts to be tight at the optimal solution of the full problem. Since all transferred cuts are valid for the full problem, they cannot exclude an optimal solution. Tightness at the surrogate incumbent is used only as a selection criterion for cuts that are likely to be useful in the full solve. A transferred cut may therefore not be tight at the final optimum while still strengthening the initial relaxation.

\subsection{When \mymethod\ helps}
\label{sec:accounting_model}
When solving $\mathcal{I}(\mathcal{B},\emptyset)$, also referred to as regular branch-and-cut, a substantial share of the runtime can be spent in the \emph{root cut loop}, separating cuts to tighten a weak initial relaxation. \mymethod\ aims to reduce this effort by injecting cuts collected from the simpler surrogate problem(s) before solving the full problem. We derive a simple model for \mymethod's speedup over regular branch-and-cut.

We model the regular branch-and-cut wall-clock time as
\begin{equation*}
    T_{\text{B\&C}} := t_{\text{root}} + t_{\text{branch}} + t_0,
\end{equation*}
where $t_{\text{root}}$ is the time spent in the root cut loop, $t_{\text{branch}}$ the time spent branching, and $t_0$ the remaining overhead. For TULIP, let $t_{\text{sur}}$ denote the time required to construct and solve the surrogate problem(s) and collect cuts. If warm-starting with these cuts changes the root cut loop time to $\alpha t_{\text{root}}$, with $\alpha\geq0$ called the root cut loop time ratio, then the TULIP wall-clock time is
\begin{equation*}
    T_{\text{TULIP}} := t_{\text{sur}} + \alpha t_{\text{root}} + t_{\text{branch}} + t_0,
\end{equation*}
where we make the simplifying assumption that the branching time and remaining overhead are unchanged, so that the model isolates \mymethod's intended effect on the root cut loop. Define $f$ as the \emph{root cut loop fraction} and $\tau$ as the
\emph{surrogate-overhead fraction}:
\begin{equation*}
    f := \frac{t_{\text{root}}}{T_{\text{B\&C}}}, \qquad \tau := \frac{t_{\text{sur}}}{T_{\text{B\&C}}},
\end{equation*}
The resulting speedup factor of \mymethod\ with respect to regular branch-and-cut is
\begin{equation}
    \eta := \frac{T_{\text{B\&C}}}{T_{\text{TULIP}}} = \frac{1}{1-(1-\alpha)f+\tau}. \label{eq:speedup-law}
\end{equation}
We refer to Equation~\eqref{eq:speedup-law} as the \mymethod\ speedup model. Equation~\eqref{eq:speedup-law} resembles Amdahl's law~\cite{amdahl1967validity}: $(1-\alpha)f$ is the fraction of the baseline runtime saved in the root cut loop, while $\tau$ is the relative cost of obtaining this saving. In particular, TULIP is beneficial compared to regular branch-and-cut when
\[
    \eta>1 \quad\Longleftrightarrow\quad (1-\alpha)f>\tau.
\]
Thus, the root cut loop time saved must exceed the surrogate overhead. Since $\tau\geq0$, a speedup requires $\alpha<1$ whenever $f>0$. Moreover, since $\alpha\geq0$ and $\tau\geq0$,
\[
    \eta = \frac{1}{1-f+\alpha f+\tau} \leq \frac{1}{1-f}, \qquad f<1.
\]
Hence, for a given baseline root cut loop fraction $f$, the speedup is bounded by
\[
    \eta_{\max} = \frac{1}{1-f},
\]
with the bound attained when $\alpha=0$ and $\tau=0$. For example, if the time spent in the root cut loop is 80\% of the total regular branch-and-cut time, the speedup factor of TULIP is at most $5$.

\Needspace{9\baselineskip}
This bound provides a baseline-only screening rule. For a target speedup $s_{\min}>1$, \mymethod\ cannot attain the target if
\[
    \frac{1}{1-f}<s_{\min}
    \quad\Longleftrightarrow\quad
    f<1-\frac{1}{s_{\min}}.
\]
Passing this test is necessary but not sufficient: deciding whether a particular \mymethod\ variant will be beneficial also requires estimates of its root cut loop time ratio $\alpha$ and surrogate-overhead fraction $\tau$.

In conclusion, the analysis identifies three conditions under which \mymethod\ is beneficial compared to regular branch-and-cut: the root cut loop time by regular branch-and-cut is relatively large compared to its total time ($f$), the time reduction factor of the root cut loop by warm-starting with TULIP's cuts is large ($1 - \alpha$), and the relative overhead of obtaining TULIP's cuts compared to regular branch-and-cut's total time is small ($\tau$).

\section{Computational study} \label{sec:computational_study}
This section describes the computational study to evaluate \mymethod\ for two-stage stochastic programs with different problem structures and exponential constraint sets. Section~\ref{sec:case_studies} introduces the case studies and describes how the generic \mymethod\ framework is applied to each, while Section~\ref{sec:experimental_design} presents the experimental design and Section~\ref{sec:performance_measures} discusses the performance measures. We end with the implementation details in Section~\ref{sec:computational_environment}. The results can be found in Section~\ref{sec:results}.

\subsection{Case studies} \label{sec:case_studies}
We evaluate four two-stage stochastic programs from the literature, see Table~\ref{tab:case-studies-methods}. In short, the \emph{stochastic Steiner tree} problem (SSTP)~\cite{leitner_decomposition_2018} connects a set of terminals in a network at minimum cost, with connectivity cuts in the second stage. The \emph{windfarm} problem, for which we build a stochastic version of the deterministic model of~\cite{pedersen2026integrated}, lays out cables connecting turbines to a substation, with first-stage connectivity cuts. The \emph{stochastic capacitated vehicle (truck) routing} problem (SCVRP)~\cite{birge_introduction_2011} serves customers with uncertain demand using capacity-limited routes, and the \emph{storage, layout and routing} problem (SLR)~\cite{narum_problem-based_2024} assigns items to warehouse locations and routes the order-picking tours under uncertain demand baskets.

For SSTP and SLR, we additionally consider variants, denoted by \sstpp\ and \slrp, that introduce an exponential number of connectivity constraints in the first stage. In contrast to SSTP, in \sstpp the first-stage edges must also connect a prescribed set of first-stage terminals. \sstpp\ is motivated by applications such as future-proof ship pipe routing, where the network installed today must already connect the current system while allowing for extensions and adaptations under future scenarios~\citep{markhorst2025future}. In \slrp, the first stage additionally installs a connected backbone linking the stocked locations, which can be used by the second-stage picking routes. For the full models, instances, and parameters of each case study, we refer to Appendices~\ref{sec:case1}-\ref{sec:slr}.

The generic \mymethod\ framework applies to all case studies, but not every surrogate construction is meaningful for every problem structure. Table~\ref{tab:case-studies-methods} summarizes the variants evaluated for each case study. For example, for \emph{windfarm}, we consider only TULIP-R as the second-stage blocks contain no cuts, while the first-stage block becomes uninformative when separated from the recourse penalty. Hence, neither
TULIP-F nor TULIP-S provides useful information, and consequently TULIP-C.

\begin{table}[ht!]
    \centering
    \caption{The case studies, their exponential constraint structure, and the \mymethod\ variants evaluated in the experiments.}
    \label{tab:case-studies-methods}
    \begin{tabular}{@{}lccccccc@{}}
        \toprule
        & \multicolumn{2}{c}{Exponential constraints}
        & \multicolumn{4}{c}{TULIP variant} 
        & \\
        \cmidrule(lr){2-3}
        \cmidrule(lr){4-7}
        Problem 
        & First stage 
        & Second stage 
        & R 
        & F 
        & S 
        & C 
        & Appendix \\
        \midrule
        SSTP     
        &            & \checkmark
        & \checkmark &            & \checkmark &
        & \ref{sec:case1} \\

        \sstpp   
        & \checkmark & \checkmark
        & \checkmark & \checkmark & \checkmark & \checkmark
        & \ref{sec:case1} \\

        Windfarm 
        & \checkmark &
        & \checkmark &            &            &
        & \ref{sec:case2} \\

        SCVRP    
        & \checkmark & \checkmark
        & \checkmark & \checkmark & \checkmark & \checkmark
        & \ref{sec:scvrp} \\

        SLR      
        &            & \checkmark
        & \checkmark &            & \checkmark &
        & \ref{sec:slr} \\

        \slrp    
        & \checkmark & \checkmark
        & \checkmark & \checkmark & \checkmark & \checkmark
        & \ref{sec:slr} \\
        \bottomrule
    \end{tabular}
\end{table}

\subsection{Experimental setup} \label{sec:experimental_design}
We use the following experimental setup. Every method has a wall-clock budget of 3600 seconds per instance. For each instance, regular branch-and-cut is run once as a standalone baseline, and each \mymethod\ variant is run independently and compared with this baseline. For \mymethod, the budget covers the complete procedure, including the surrogate phase and the final full-problem solve. Whenever the surrogate solve(s) yield mutually consistent values for decision variables of the full problem, these values are provided as a primal start for the final solve; otherwise, only the collected cuts are transferred.

For TULIP-R, we use Feng and Ryan's problem-based scenario reduction method~\citep{feng2016solution}. We choose this method because it is computationally inexpensive and provides a practical compromise between reduction quality and runtime, which is important for \mymethod\ because the scenario-reduction time directly contributes to the surrogate overhead. The probabilities $\tilde p^{(s)}$ of the retained scenarios are those obtained from the scenario-reduction procedure. To account for the method's stochasticity, we run TULIP-R with three different seeds on each instance. The reduced scenario set size is $|\mathcal{S}'| = \max(2,\lceil 0.1\,|\mathcal{S}| \rceil)$ for all case studies except windfarm; see Appendix~\ref{sec:case2}. The baseline and the deterministic TULIP-F, TULIP-S, and TULIP-C are run once. For \sstpp, we also report SAA results using the same $10\%$ scenario-set size as a simple benchmark. This provides context on the solution quality attainable from a random scenario subset without exploiting problem-specific information. We perform three independent SAA replications per instance and evaluate each resulting first-stage solution on the full scenario set; further details are provided in Appendix~\ref{sec:case1_results}. Case-study-specific experimental settings are given in Appendices~\ref{sec:experiments1}, \ref{sec:experiments2}, \ref{sec:experiments3}, and~\ref{sec:experiments4} for SSTP/\sstpp, windfarm, SCVRP, and SLR/\slrp, respectively.

\subsection{Performance measures} \label{sec:performance_measures}

Throughout the computational study, we use \emph{performance} as an umbrella term for the three measures considered below: the ability to solve instances to optimality, runtime performance on instances solved by all methods, and final optimality gaps on instances solved by no method. 

For ease of presentation, let $\mathcal{J}$ denote the set of
instances considered for each group of experiments and let $\mathcal{M}$ denote the set of methods
compared. The latter consists of regular branch-and-cut and the TULIP variants
applicable to the corresponding problem, as indicated in
Table~\ref{tab:case-studies-methods}. For a method $m \in \mathcal{M}$,
let
\[
    \mathcal{J}_m^\star := \{j \in \mathcal{J} : m \text{ solves instance } j \text{ to optimality within the time limit}\}.
\]
We consider an instance solved to optimality when its final optimality
gap is at most $0.01\%$, as this is the default in the Gurobi solver. We further define
\[
    \mathcal{J}^{\cap} := \bigcap_{m \in \mathcal{M}} \mathcal{J}_m^\star,
    \qquad
    \mathcal{J}^{0} := \mathcal{J} \setminus \bigcup_{m \in \mathcal{M}} \mathcal{J}_m^\star,
\]
as the sets of instances solved by all compared methods and by no
compared method, respectively.

Since our main focus is on solving instances to optimality, we first report the number and percentage of instances solved to optimality by each method, where the percentage is computed as $|\mathcal{J}_m^\star|/|\mathcal{J}| \cdot 100$. Second, for the instances solved to optimality by all methods, we compare how quickly optimality is reached. For a fair runtime comparison, we restrict this comparison to $\mathcal{J}^{\cap}$, such that all methods are evaluated on exactly the same instances. The generic runtimes from Section~\ref{sec:accounting_model} are indexed here by instance and method, with $T_{j,m}$ denoting the observed wall-clock times for method $m$ on instance $j$. For an instance $j \in \mathcal{J}^{\cap}$ and method $m$, we define the observed speedup regular over branch-and-cut as
\[
    \eta^{\mathrm{obs}}_{j,m} := \frac{T_{j,\mathrm{B\&C}}}{T_{j,m}}.
\]
Hence, $\eta^{\mathrm{obs}}_{j,m} > 1$ indicates that method $m$ is faster than regular branch-and-cut. For reporting results over multiple instances, we aggregate these quantities over $\mathcal{J}^{\cap}$. For method $m \in \mathcal{M}$, we define the mean runtime as
\[
    \bar{T}_m := \frac{1}{|\mathcal{J}^{\cap}|} \sum_{j \in \mathcal{J}^{\cap}} T_{j,m},
\]
and, for each TULIP variant $m$, the mean observed speedup as
\[
    \bar{\eta}_m := \frac{1}{|\mathcal{J}^{\cap}|} \sum_{j \in \mathcal{J}^{\cap}} \eta^{\mathrm{obs}}_{j,m}.
\]
Finally, for the instances in $\mathcal{J}^{0}$ that are not solved to optimality by any method, a runtime comparison is not meaningful. We therefore compare the final optimality gaps as a proxy for how much
progress each method makes toward optimality. For method $m$, we report the mean final gap
\[
    \bar{g}_m := \frac{1}{|\mathcal{J}^{0}|} \sum_{j \in \mathcal{J}^{0}} g_{j,m},
\]
where $g_{j,m}$ denotes the final gap for instance $j$ by method $m$. We additionally define the instance-level gap difference relative to regular branch-and-cut as
\[
    \Delta g_{j,m} := g_{j,\mathrm{B\&C}} - g_{j,m},
\]
so that $\Delta g_{j,m}>0$ indicates that method $m$ terminates with a smaller optimality gap than regular branch-and-cut. Over the instances in $\mathcal{J}^0$, we report the mean gap difference
\[
    \overline{\Delta g}_m := \frac{1}{|\mathcal{J}^0|} \sum_{j\in\mathcal{J}^0} \Delta g_{j,m}.
\]

For TULIP-R, which is randomized through the scenario-reduction step, we perform three runs with different seeds. Let $\mathcal{R}_{\mathrm{seed}}=\{1,2,3\}$ denote these runs, and let $T_{j,\mathrm{TULIP\text{-}R}}^{(r)}$ and $g_{j,\mathrm{TULIP\text{-}R}}^{(r)}$ denote the runtime and final gap of run $r$ on instance $j$, respectively. We aggregate these runs into the instance-level quantities used above as
\[
    T_{j,\mathrm{TULIP\text{-}R}} := \frac{1}{|\mathcal{R}_{\mathrm{seed}}|}\sum_{r\in\mathcal{R}_{\mathrm{seed}}} T_{j,\mathrm{TULIP\text{-}R}}^{(r)},
    \qquad
    g_{j,\mathrm{TULIP\text{-}R}} := \frac{1}{|\mathcal{R}_{\mathrm{seed}}|}\sum_{r\in\mathcal{R}_{\mathrm{seed}}}
    g_{j,\mathrm{TULIP\text{-}R}}^{(r)}.
\]
An instance is considered solved by TULIP-R only if all three runs attain a final gap of at most $0.01\%$. For reporting purposes, gaps are capped at $100\%$ to limit the influence of extreme values on the averages. A run without an incumbent is assigned the capped value of $100\%$. All subsequent performance measures are computed from these instance-level quantities, with TULIP-R contributing one observation per instance, as do the deterministic methods.

\subsection{Implementation details} \label{sec:computational_environment}
All experiments were run using Julia~1.10.4 and Gurobi~10.0.1 via \texttt{Stochastic\allowbreak Programs.jl}~\citep{spjl} and \texttt{JuMP}~\citep{Lubin2023} on the Snellius national supercomputer operated by SURF, on compute nodes equipped with AMD EPYC 9654 processors. Each job was allocated $24$\,GB of memory, except for windfarm, which required $64$\,GB to avoid running out of memory on the larger instances. In total, the experiments comprise $3{,}076$ runs: SSTP ($949$), SLR ($1{,}522$), SCVRP ($270$), windfarm ($159$), SAA ($154$), and EVPI ($22$). Gurobi is configured as single-threaded for all runs to isolate the algorithmic effect of TULIP from effects of parallel solver execution. See our GitHub repository for other implementation details~\cite{markhorst_tulipjl}.

\section{Results} \label{sec:results}
Table~\ref{tab:summary-all} summarizes the results for each case study, pooling all instances across the tested scenario-set sizes. It reports $|\mathcal{J}_m^\star|$, the number of instances solved to optimality by method $m$; $\bar T_m$ and $\bar\eta_m$, the mean runtime and speedup with respect to regular branch-and-cut on $\mathcal{J}^\cap$; and $\bar g_m$, the mean final optimality gap on $\mathcal{J}^0$. More detailed results by scenario count are provided in Appendices~\ref{sec:case1}--\ref{sec:slr}.

\begin{table}[ht!]
\centering
\caption{Overview of the numerical results. Bold values indicate the best performance within each problem and performance measure, and ties are all highlighted. The speedup of branch-and-cut relative to itself is $1.00$ by definition.}
\label{tab:summary-all}
\small
\setlength{\tabcolsep}{5pt}
\begin{tabular}{llrcrrcrc}
\toprule
&&& & \multicolumn{3}{c}{Solved by all methods} & \multicolumn{2}{c}{Solved by no method} \\
\cmidrule(lr){5-7}\cmidrule(lr){8-9}
Problem & Method $m$ & $|\mathcal{J}|$ & $|\mathcal{J}_m^\star|$ & $|\mathcal{J}^{\cap}|$ & $\bar{T}_m$ (s) & $\bar{\eta}_m$ & $|\mathcal{J}^{0}|$ & $\bar{g}_m$ (\%) \\
\midrule
SSTP+ & B\&C & 77 & 72\,(94\%) & 71 & 245 & 1.00 & 2 & 100.0 \\
 & TULIP-R &  & 74\,(96\%) &  & 149 & 1.10 &  & 100.0 \\
 & TULIP-F &  & \textbf{75\,(97\%)} &  & \textbf{122} & \textbf{1.73} &  & 100.0 \\
 & TULIP-S &  & 72\,(94\%) &  & 309 & 0.76 &  & 100.0 \\
 & TULIP-C &  & 72\,(94\%) &  & 214 & 1.15 &  & \textbf{54.4} \\
\addlinespace

SSTP & B\&C & 77 & \textbf{77\,(100\%)} & 77 & 102 & 1.00 & 0 & -- \\
 & TULIP-R &  & \textbf{77\,(100\%)} &  & 122 & 0.60 &  & -- \\
 & TULIP-S &  & \textbf{77\,(100\%)} &  & \textbf{83} & \textbf{1.07} &  & -- \\
\addlinespace

SLR & B\&C & 80 & \textbf{67\,(84\%)} & 62 & \textbf{149} & \textbf{1.00} & 13 & 16.0 \\
 & TULIP-R &  & 62\,(78\%) &  & 153 & 0.66 &  & \textbf{12.4} \\
 & TULIP-S &  & \textbf{67\,(84\%)} &  & 153 & 0.88 &  & 16.1 \\
\addlinespace

SLR+ & B\&C & 60 & \textbf{50\,(83\%)} & 48 & \textbf{508} & \textbf{1.00} & 10 & \textbf{12.6} \\
 & TULIP-R &  & 48\,(80\%) &  & 512 & 0.94 &  & 20.4 \\
 & TULIP-F &  & \textbf{50\,(83\%)} &  & 510 & 0.98 &  & \textbf{12.6} \\
 & TULIP-S &  & \textbf{50\,(83\%)} &  & 513 & 0.97 &  & 12.7 \\
 & TULIP-C &  & \textbf{50\,(83\%)} &  & 515 & 0.96 &  & 12.7 \\
\addlinespace

Windfarm & B\&C & 18 & 12\,(67\%) & 12 & \textbf{336} & \textbf{1.00} & 5 & 60.2 \\
 & TULIP-R &  & \textbf{13\,(72\%)} &  & 447 & 0.69 &  & \textbf{56.4} \\
\addlinespace

SCVRP & B\&C & 27 & 10\,(37\%) & 10 & \textbf{65} & \textbf{1.00} & 16 & 11.3 \\
 & TULIP-R &  & 10\,(37\%) &  & 77 & 0.18 &  & 7.2 \\
 & TULIP-F &  & 10\,(37\%) &  & 228 & 0.43 &  & \textbf{4.5} \\
 & TULIP-S &  & \textbf{11\,(41\%)} &  & \textbf{65} & 0.78 &  & 11.3 \\
 & TULIP-C &  & 10\,(37\%) &  & 273 & 0.38 &  & \textbf{4.5} \\
\bottomrule
\end{tabular}
\end{table}
Overall, considering $\mathcal{J}_m^\star$, the \mymethod\ variants solve at least as many instances to optimality within the time limit as regular branch-and-cut. Some variants solve slightly more instances for \sstpp, windfarm,
and SCVRP. The exceptions are TULIP-R on SLR and \slrp, where
fewer instances are solved. Since TULIP-R is, as described in Section~\ref{sec:experimental_design}, randomized through the scenario-reduction step, we examine the sensitivity of its performance to the scenario-reduction seed separately in Appendix~\ref{sec:stability_tulip_r}.

The runtime comparison on $\mathcal{J}^\cap$ shows a clearer distinction between the case studies. The strongest gains occur for \sstpp\, where TULIP-F is substantially faster than regular branch-and-cut, and TULIP-R and TULIP-C also achieve a mean speedup above one. For SSTP, only TULIP-S gives a small speedup, whereas TULIP-R is slower. In the remaining case studies, the \mymethod\ variants are generally slower than regular branch-and-cut or close to break-even. Thus, we see the main runtime benefit of \mymethod\ in the stochastic Steiner tree problems, and particularly in \sstpp. This comparison should, however, be interpreted in light of the experimental design. To allow for an extensive computational comparison, we selected relatively tractable instances: among those solved to optimality, the average runtime was below $10$ minutes despite the $1$-hour time limit. For such instances, the surrogate overhead introduced by \mymethod\ is relatively large and therefore harder to compensate for. Further inspection shows that \mymethod\ is mainly slower on the smaller SLR instances, which lowers the mean speedup factor. On the more difficult instances, the difference becomes much smaller, and \mymethod\ occasionally recovers time. Since these instances dominate the total runtime, the mean runtimes of \mymethod\ and regular branch-and-cut remain similar.

Finally, on $\mathcal{J}^0$, runtime comparisons are not meaningful because no method reaches optimality within the time limit. The final optimality gaps nevertheless show that some \mymethod\ variants make more progress toward optimality than regular branch-and-cut. Several \mymethod\ variants terminate with smaller final optimality gaps than branch-and-cut, most notably TULIP-C for \sstpp\ and TULIP-R, TULIP-F, and TULIP-C for SCVRP, while for SLR, \slrp, and windfarm, the differences are smaller or mixed.

The remainder of this section examines these differences in more detail. In Section~\ref{sec:where}, we first discuss where \mymethod\ improves upon regular branch-and-cut, contrasting the stochastic Steiner tree case studies with the remaining problems. In Section~\ref{sec:why}, we then explain the observed runtime differences using the speedup model from Section~\ref{sec:accounting_model}, focusing on the root cut loop time saved and the surrogate overhead. We also use these results to address our third research question on the choice of surrogate.

\subsection{Where \mymethod\ improves upon regular branch-and-cut} \label{sec:where}
The results show that \mymethod\ can considerably improve upon regular branch-and-cut for some problems, while providing little or no benefit for others. We also examine how this comparison changes with the number of scenarios $|\mathcal{S}|$. This is relevant because increasing $|\mathcal{S}|$ enlarges the full stochastic problem, while the computational growth of the \mymethod\ surrogate depends on its construction. Hence, the relative benefit of the warm start may vary with the scenario-set size. We discuss these differences across case studies and scenario-set sizes separately below. We therefore examine these patterns across problems and scenario-set sizes below. For conciseness, we report the detailed results for \sstpp\ in the main text. The corresponding detailed results for SSTP, windfarm, SCVRP, and SLR/\slrp\ are provided in Appendices~\ref{sec:case1_results}, \ref{sec:case2_results}, \ref{sec:case3_results}, and \ref{sec:case4_results}, respectively.

\subsubsection{Strong improvements: the stochastic Steiner tree problem}
Table~\ref{tab:solved-speedup-steiner-full} first shows that the methods solve a similar number of \sstpp\ instances to optimality. TULIP-F performs best in this respect, with $|\mathcal{J}_{\mathrm{F}}^\star|$ exceeding $|\mathcal{J}_{\mathrm{B\&C}}^\star|$ at $|\mathcal{S}|=150$ and $300$, while the differences are small for the remaining scenario counts. On the commonly solved instances $\mathcal{J}^{\cap}$, TULIP-F consistently outperforms regular branch-and-cut in runtime, with $\bar{\eta}_{\mathrm{F}}$ increasing from $1.32$ at $|\mathcal{S}|=50$ to $2.56$ at $|\mathcal{S}|=500$. TULIP-R and TULIP-C also become increasingly effective as the number of scenarios grows, reaching $\bar{\eta}_{\mathrm{R}}=1.72$ and $\bar{\eta}_{\mathrm{C}}=1.89$, respectively, at $|\mathcal{S}|=500$. In contrast, $\bar{\eta}_{\mathrm{S}}<1$ for all scenario counts. Figure~\ref{fig:steiner-nodes-full} in Appendix~\ref{sec:case1_results} shows a general tendency toward larger speedups at higher scenario counts, although this pattern is not consistent for every graph size.

\begin{table}[ht]
\centering
\caption{Number of \sstpp\ instances solved and mean speedup by scenario count. Bold values indicate the best performance within each scenario count; ties are all highlighted.}
\label{tab:solved-speedup-steiner-full}
\footnotesize
\setlength{\tabcolsep}{2.8pt}
\begin{tabular}{rrrcccccccccc}
\toprule
&&& \multicolumn{2}{c}{B\&C} & \multicolumn{2}{c}{TULIP-R} & \multicolumn{2}{c}{TULIP-F} & \multicolumn{2}{c}{TULIP-S} & \multicolumn{2}{c}{TULIP-C} \\
\cmidrule(lr){4-5}\cmidrule(lr){6-7}\cmidrule(lr){8-9}\cmidrule(lr){10-11}\cmidrule(lr){12-13} 
$|\mathcal{S}|$ & $|\mathcal{J}|$ & $|\mathcal{J}^{\cap}|$ & $|\mathcal{J}_{\mathrm{B\&C}}^\star|$ & $\bar{T}_{\mathrm{B\&C}}$ (s) & $|\mathcal{J}_\mathrm{R}^\star|$ & $\bar{\eta}_{\mathrm{R}}$ & $|\mathcal{J}_\mathrm{F}^\star|$ & $\bar{\eta}_\mathrm{F}$ & $|\mathcal{J}_\mathrm{S}^\star|$ & $\bar{\eta}_\mathrm{S}$ & $|\mathcal{J}_\mathrm{C}^\star|$ & $\bar{\eta}_\mathrm{C}$ \\
\midrule
$50$ & 11 & 11 &
\textbf{11\,(100\%)} & 38 &
\textbf{11\,(100\%)} & 0.62 &
\textbf{11\,(100\%)} & \textbf{1.32} &
\textbf{11\,(100\%)} & 0.76 &
\textbf{11\,(100\%)} & 0.92 \\

$75$ & 11 & 11 &
\textbf{11\,(100\%)} & 99 &
\textbf{11\,(100\%)} & 0.79 &
\textbf{11\,(100\%)} & \textbf{1.46} &
\textbf{11\,(100\%)} & 0.80 &
\textbf{11\,(100\%)} & 1.02 \\

$100$ & 11 & 11 &
\textbf{11\,(100\%)} & 150 &
\textbf{11\,(100\%)} & 0.93 &
\textbf{11\,(100\%)} & \textbf{1.60} &
\textbf{11\,(100\%)} & 0.69 &
\textbf{11\,(100\%)} & 1.03 \\

$150$ & 11 & 10 &
10\,(91\%) & 174 &
\textbf{11\,(100\%)} & 1.09 &
\textbf{11\,(100\%)} & \textbf{1.73} &
\textbf{11\,(100\%)} & 0.75 &
10\,(91\%) & 1.11 \\

$200$ & 11 & 10 &
\textbf{11\,(100\%)} & 313 &
\textbf{11\,(100\%)} & 1.26 &
\textbf{11\,(100\%)} & \textbf{1.72} &
10\,(91\%) & 0.70 &
10\,(91\%) & 1.03 \\

$300$ & 11 & 9 &
9\,(82\%) & 200 &
10\,(91\%) & 1.45 &
\textbf{11\,(100\%)} & \textbf{1.86} &
9\,(82\%) & 0.76 &
10\,(91\%) & 1.17 \\

$500$ & 11 & 9 &
\textbf{9\,(82\%)} & 838 &
\textbf{9\,(82\%)} & 1.72 &
\textbf{9\,(82\%)} & \textbf{2.56} &
\textbf{9\,(82\%)} & 0.85 &
\textbf{9\,(82\%)} & 1.89 \\
\bottomrule
\end{tabular}
\end{table}

Table~\ref{tab:gapdiff-steiner-full} considers the instances in $\mathcal{J}^{0}$, i.e., those not solved to optimality by any method within the time limit. Such instances occur only at $|\mathcal{S}|=500$, where $|\mathcal{J}^{0}|=2$ and regular branch-and-cut terminates without an incumbent, yielding $\bar{g}_{\mathrm{B\&C}}=100\%$. TULIP-R, TULIP-F, and TULIP-S do not improve on this final gap, with $\overline{\Delta g}_{\mathrm{R}}=\overline{\Delta g}_{\mathrm{F}}=\overline{\Delta g}_{\mathrm{S}}=0.0$, whereas TULIP-C achieves $\overline{\Delta g}_{\mathrm{C}}=45.6$ percentage points. Thus, although TULIP-F gives the greatest runtime improvements on $\mathcal{J}^{\cap}$, TULIP-C makes the most progress on the hardest instances in $\mathcal{J}^{0}$.

\begin{table}[ht]
\centering
\caption{Final optimality gaps for \sstpp\ instances not solved by any method within the time limit. Bold values indicate the best performance within each scenario count.}
\label{tab:gapdiff-steiner-full}
\small
\begin{tabular}{crccccc}
\toprule
$|\mathcal{S}|$ & $|\mathcal{J}^0|$ & $\bar{g}_{\mathrm{B\&C}}$ (\%) & $\overline{\Delta g}_\mathrm{R}$ & $\overline{\Delta g}_\mathrm{F}$ & $\overline{\Delta g}_\mathrm{S}$ & $\overline{\Delta g}_\mathrm{C}$ \\
\midrule
$50$-$300$ & 0 & -- & -- & -- & -- & -- \\
$500$ & 2 & 100.0 & +0.0 & +0.0 & +0.0 & \textbf{+45.6} \\
\bottomrule
\end{tabular}
\end{table}

\begin{figure}[ht]
    \centering
    \includegraphics[width=0.75\linewidth]{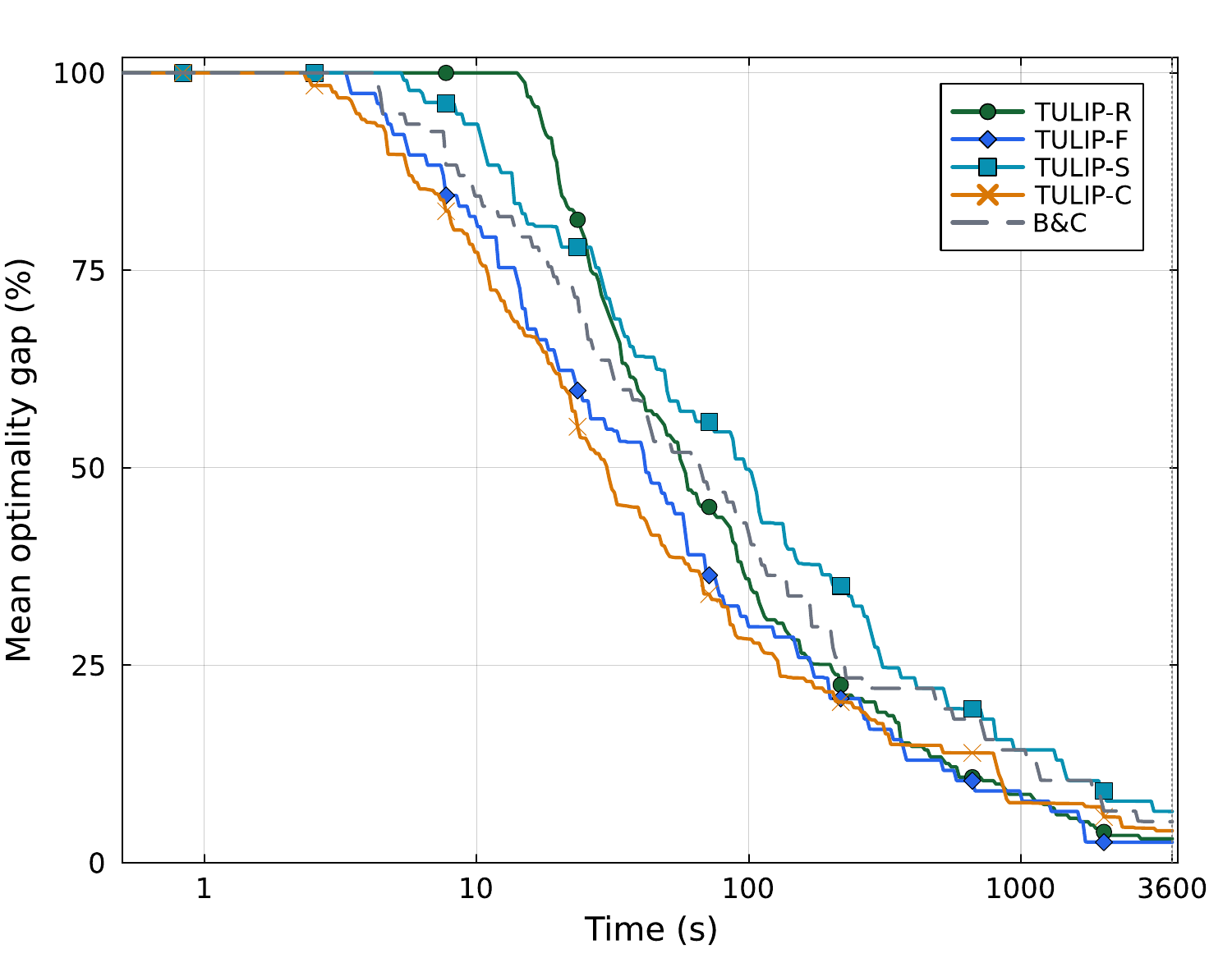}
    \caption{Mean optimality gap over time for \sstpp.}
    \label{fig:steiner_full}
\end{figure}

Figure~\ref{fig:steiner_full} complements these performance measures by showing how the mean optimality gap develops over time across all \sstpp\ instances. The figure reflects the same overall pattern: TULIP-F reduces the gap more quickly than regular branch-and-cut over most of the solve, whereas TULIP-S generally makes slower progress; TULIP-R and TULIP-C lie between these two extremes. Hence, the benefit of TULIP-F is not limited to its larger $\bar{\eta}_{\mathrm{F}}$ on $\mathcal{J}^{\cap}$ but is visible throughout the optimization process. For completeness, Appendix~\ref{sec:case1_results} also reports the \textit{expected value of perfect information} (EVPI) for the \sstpp\ instances: at $|\mathcal{S}|=50$, the mean instance-level EVPI percentage is $6.8\%$, confirming that the uncertainty has a non-negligible effect on the optimal solution.

The results for SSTP are more modest. Table~\ref{tab:solved-speedup-steiner-nofs} in Appendix~\ref{sec:case1_results} shows that all $77$ instances are solved to optimality by regular branch-and-cut, TULIP-R, and TULIP-S, so $\mathcal{J}^0 = \emptyset$. TULIP-S performs close to regular branch-and-cut for smaller scenario sets and becomes increasingly beneficial as $|\mathcal{S}|$ grows, with its mean speedup increasing from $0.99$ at $|\mathcal{S}|=50$ to $1.18$ at $|\mathcal{S}|=500$. TULIP-R, in contrast, remains slower than regular branch-and-cut for every scenario count, although its relative performance likewise improves from $0.37$ to $0.83$. Figure~\ref{fig:steiner_nofs} in Appendix~\ref{sec:case1_results} confirms this pattern over the course of the solve: TULIP-S follows regular branch-and-cut closely, whereas TULIP-R generally reduces the optimality gap more slowly. Thus, unlike for \sstpp, the warm start yields only a modest speedup for SSTP, and only when cuts are collected from the individual second-stage blocks through TULIP-S.

\subsubsection{Other case studies: limited improvements} \label{subsubsec:other_cases}
For windfarm, TULIP-R does not provide a systematic speedup over regular branch-and-cut. Although it solves slightly more instances overall, its mean speedup on the instances solved by both methods remains below one. In terms of the parameters introduced in Section~\ref{sec:accounting_model}, this indicates that the reduction in the root cut loop fraction, represented by $(1-\alpha)f$, is not large enough to compensate for the surrogate-overhead fraction $\tau$. In particular, the measured root cut loop fraction $f$ is small, and the direct root-time saving $(1-\alpha)f$ is slightly negative; these quantities are reported in Section~\ref{sec:why} and Table~\ref{tab:root_bottleneck}. Their effectiveness also depends strongly on the first-stage turbine layout. Consequently, the root cut loop savings $(1-\alpha)f$ are insufficient to offset the surrogate-overhead fraction $\tau$. Detailed results are provided in Appendix~\ref{sec:case2_results}.

The SCVRP provides an even clearer example of the importance of the root cut loop fraction. None of the \mymethod\ variants achieves a mean speedup greater than $1$, although some obtain smaller final gaps on instances that remain unsolved. The transferred first-stage subtour cuts can strengthen the initial relaxation, but regular branch-and-cut already closes the root considerably quickly. Thus, the root cut loop savings $(1-\alpha)f$ are too small to offset the surrogate-overhead fraction $\tau$. Detailed results are reported in Section~\ref{sec:case3_results}.

For SLR and \slrp, the \mymethod\ variants also remain close to or below the performance of regular branch-and-cut, with only occasional small speedups for individual scenario counts. The main limitation is that the second-stage subtour cuts depend strongly on the first-stage assignment. As a result, the first logged root-bounds show little improvement of the full problem's initial relaxation, so $\alpha$ remains close to $1$ and the resulting root cut loop savings $(1-\alpha)f$ are too small to offset the surrogate-overhead fraction $\tau$. The detailed results in Appendix~\ref{sec:case4_results} show that the effect of \mymethod\ on the final optimality gaps is also mixed.

In short, \mymethod\ accelerates branch-and-cut when the root cut loop time saved by the warm start is sufficiently large to outweigh the cost of constructing the surrogate. In our experiments, this occurs only for the stochastic Steiner tree problems, and most clearly for \sstpp. Section~\ref{sec:why} examines whether the \mymethod\ speedup model from Section~\ref{sec:accounting_model} quantitatively explains this pattern.

\FloatBarrier
\subsection{Why \mymethod\ works} \label{sec:why}
We now test whether Equation~\eqref{eq:speedup-law} explains the observed runtime performance of \mymethod\ across the case studies. According to the equation, TULIP is beneficial when the root cut loop time saved by the warm start exceeds the surrogate overhead, i.e.,
\[
    (1-\alpha)f > \tau.
\]
We therefore examine whether the observed speedups are consistent with the three quantities identified in Section~\ref{sec:accounting_model}: the baseline root cut loop fraction $f$, the reduction in root cut loop time $1-\alpha$, and the surrogate-overhead fraction $\tau$.

For each instance $j$, let $f_j$ denote the root cut loop fraction of the regular branch-and-cut baseline. Since the same regular branch-and-cut baseline is used for all \mymethod\ variants, $f_j$ does not depend on the variant $m$. For \mymethod\ variant $m$, let $\alpha_{j,m}$ and $\tau_{j,m}$ denote the root cut loop time ratio and surrogate-overhead fraction, respectively. Equation~\eqref{eq:speedup-law} then gives the predicted speedup
\[
    \eta^{\mathrm{pred}}_{j,m} := \frac{1}{1-(1-\alpha_{j,m})f_j+\tau_{j,m}}.
\]
We compare this quantity with the observed speedup $\eta^{\mathrm{obs}}_{j,m}$ defined in Section~\ref{sec:performance_measures}.

\begin{figure}[!htbp]
    \centering
    \begin{subfigure}{0.80\linewidth}
        \centering
        \includegraphics[width=\linewidth,trim=24bp 18bp 47bp 12bp,clip]{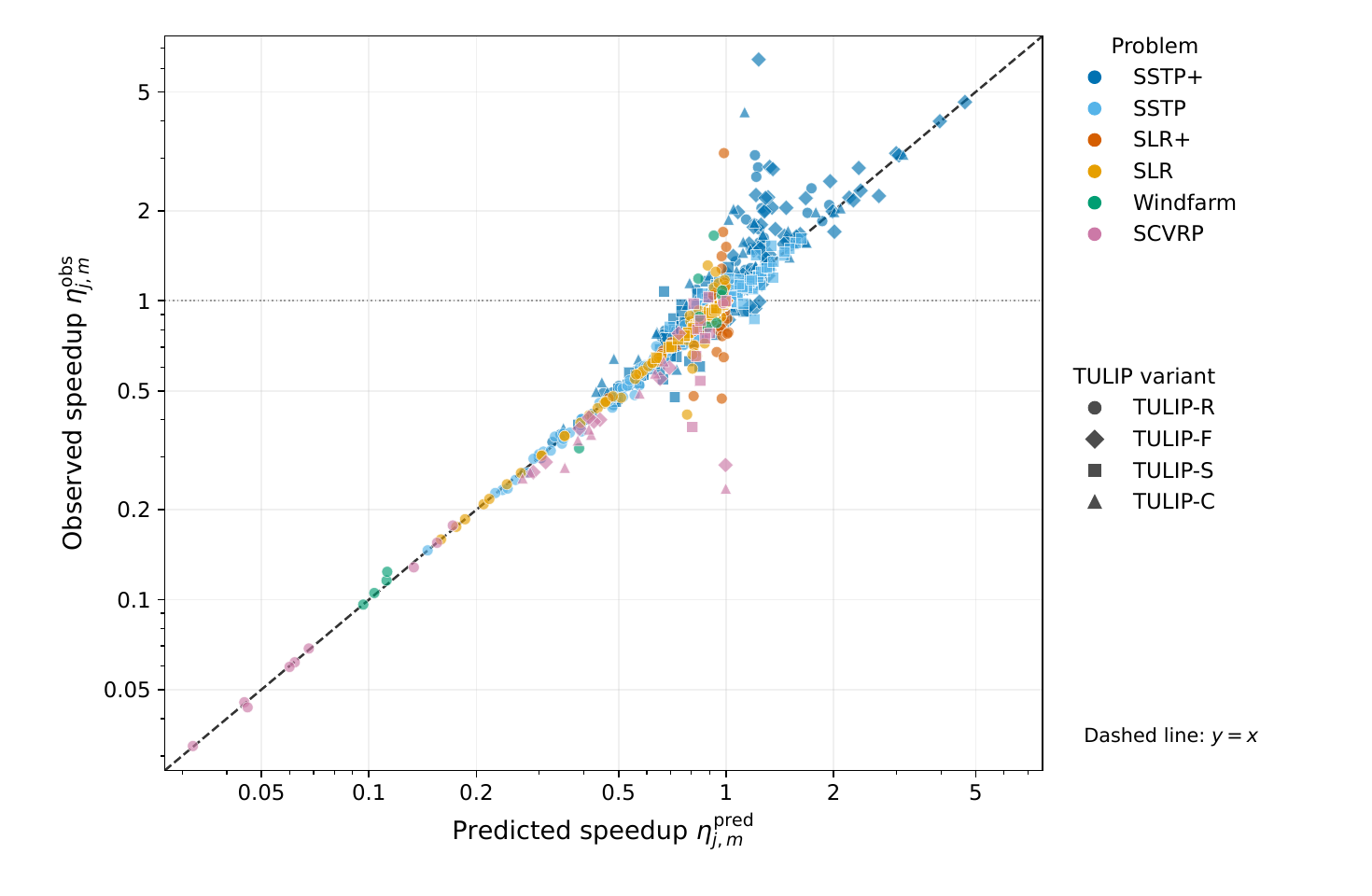}
       \subcaption{Predicted speedup from Equation~\eqref{eq:speedup-law} versus observed speedup. The dashed line indicates perfect agreement.}
        \label{fig:figC}
    \end{subfigure}

    \smallskip

    \begin{subfigure}{\linewidth}
        \centering
        \includegraphics[width=\linewidth,trim=16bp 21bp 63bp 23bp,clip]{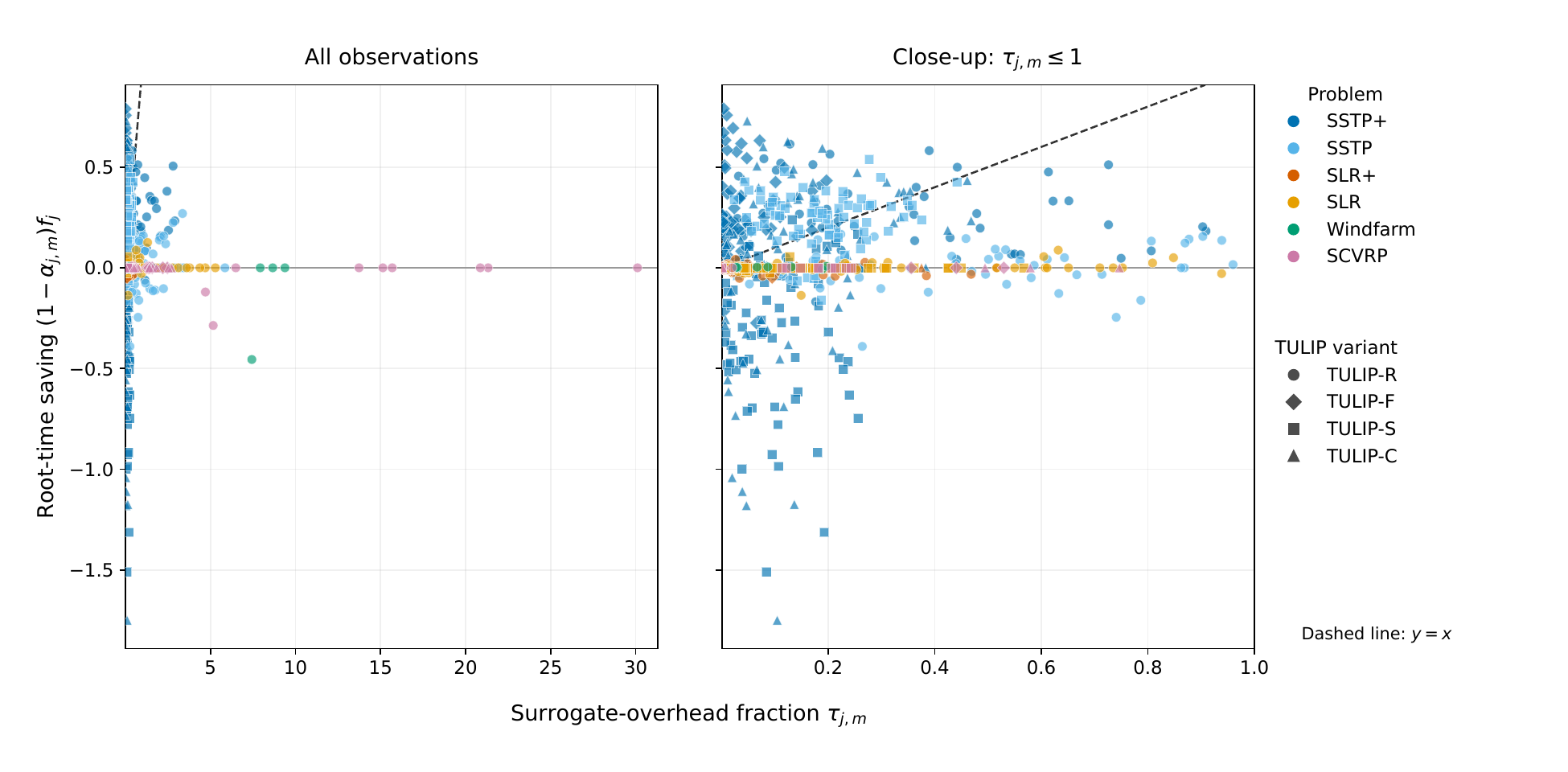}
        \subcaption{Root cut loop time saving $(1-\alpha_{j,m})f_j$ versus surrogate-overhead fraction $\tau_{j,m}$. Points above the dashed line satisfy $(1-\alpha_{j,m})f_j>\tau_{j,m}$ and are therefore predicted by Equation~\eqref{eq:speedup-law} to yield a speedup.}
        \label{fig:figD}
    \end{subfigure}

    \caption{Evaluation of the \mymethod\ speedup model across all case studies on $\mathcal{J}^{\cap}$. Figure~\ref{fig:figC} compares predicted and observed speedups, and Figure~\ref{fig:figD} shows the break-even condition $(1-\alpha_{j,m})f_j>\tau_{j,m}$.}
    \label{fig:result_c_and_d}
\end{figure}

Figure~\ref{fig:result_c_and_d} shows that the speedup model captures the main variation in \mymethod's observed speedup across the case studies. Both graphs use the instance--variant pairs from $\mathcal{J}^{\cap}$ in Table~\ref{tab:summary-all}. Figure~\ref{fig:figC} compares the speedup predicted by Equation~\eqref{eq:speedup-law} with the observed speedup. Predicted and observed speedups are strongly positively correlated ($r=0.82$, pooled over $806$ pairs), with most observations also lying reasonably close to the diagonal. Figure~\ref{fig:figD} then visualizes the break-even condition for the same observations. By construction, observations above the dashed line satisfy $(1-\alpha_{j,m})f_j>\tau_{j,m}$ and are therefore predicted to yield a speedup, whereas observations below the line are predicted not to. The model correctly predicts whether a speedup occurs for 731 of the 806 observations (90.7\%). Hence, Equation~\eqref{eq:speedup-law} provides a good explanation of both the magnitude of \mymethod's performance and its improvement over regular branch-and-cut. Additional supporting analyses are provided in Appendix~\ref{sec:figures_appendix}: Figure~\ref{fig:figA} shows that larger improvements in the first root bound reported by Gurobi tend to coincide with larger observed speedups, while Figure~\ref{fig:figB} shows the direct relationship between root cut loop speedup and overall speedup.

\begin{table}[!htb]
    \centering
    \caption{Diagnostic components of the \mymethod\ speedup model by case study and \mymethod\ variant on $\mathcal{J}^{\cap}$.}
    \label{tab:root_bottleneck}
    \begin{tabular}{llrrrrr}
\toprule
Problem & Method & $\bar f$ & $\overline{(1-\alpha)f}$ & $\bar\tau$ & $\bar\eta^{\mathrm{obs}}$ & $\bar\eta^{\mathrm{pred}}$ \\
\midrule
SSTP\textsuperscript{+} & TULIP-R & 0.788 & 0.260 & 0.589 & 1.095 & 0.935 \\
SSTP\textsuperscript{+} & TULIP-F & 0.788 & 0.267 & 0.038 & 1.725 & 1.463 \\
SSTP\textsuperscript{+} & TULIP-S & 0.788 & -0.325 & 0.101 & 0.757 & 0.740 \\
SSTP\textsuperscript{+} & TULIP-C & 0.788 & -0.006 & 0.140 & 1.151 & 1.017 \\
SSTP & TULIP-R & 0.869 & 0.014 & 0.988 & 0.600 & 0.602 \\
SSTP & TULIP-S & 0.869 & 0.231 & 0.191 & 1.074 & 1.074 \\
SLR\textsuperscript{+} & TULIP-R & 0.163 & -0.005 & 0.113 & 0.943 & 0.913 \\
SLR\textsuperscript{+} & TULIP-F & 0.163 & -0.001 & 0.022 & 0.977 & 0.979 \\
SLR\textsuperscript{+} & TULIP-S & 0.163 & -0.000 & 0.028 & 0.974 & 0.974 \\
SLR\textsuperscript{+} & TULIP-C & 0.163 & -0.000 & 0.048 & 0.955 & 0.957 \\
SLR & TULIP-R & 0.416 & 0.003 & 0.944 & 0.661 & 0.657 \\
SLR & TULIP-S & 0.416 & 0.000 & 0.150 & 0.880 & 0.879 \\
Windfarm & TULIP-R & 0.111 & -0.037 & 2.974 & 0.690 & 0.599 \\
SCVRP & TULIP-R & 0.133 & -0.041 & 13.320 & 0.175 & 0.175 \\
SCVRP & TULIP-F & 0.133 & -0.000 & 1.158 & 0.434 & 0.536 \\
SCVRP & TULIP-S & 0.133 & 0.000 & 0.164 & 0.777 & 0.862 \\
SCVRP & TULIP-C & 0.133 & 0.000 & 1.330 & 0.379 & 0.500 \\
\bottomrule
\end{tabular}

\end{table}

Table~\ref{tab:root_bottleneck} shows how the same quantities explain both the positive and negative cases. For \sstpp, TULIP-F combines a large baseline root cut loop fraction ($\bar f=0.788$) with a substantial mean root-time saving ($\overline{(1-\alpha)f}=0.267$) and negligible surrogate overhead ($\bar\tau=0.038$), resulting in the largest observed speedup. TULIP-R achieves a similar root-time saving ($0.260$), but at much higher overhead, while TULIP-S increases root time substantially on average
($\overline{(1-\alpha)f}=-0.325$), whereas the increase for
TULIP-C is marginal ($-0.006$). For SSTP, the large baseline root cut loop fraction likewise creates substantial room for improvement, but only TULIP-S exploits it: its mean root-time saving is $0.231$, whereas TULIP-R saves only $0.014$ and incurs considerable overhead. The remaining case studies fail for different reasons. \slrp and windfarm have relatively small baseline root cut loop fractions, leaving little runtime to save. For windfarm, this limitation is even worse due to substantial surrogate overhead and a slightly negative mean
root-time saving, while for SLR the mean root-time saving is approximately
zero. In SCVRP, the mean root-time saving is likewise zero or negative, and the surrogate overhead can be substantial. Thus, the negative cases are not explained by a single structural deficiency, but by different ways in which the condition $(1-\alpha)f>\tau$ fails.

The speedup model is not exact, however. Although predicted and observed speedups generally align, some observations deviate from the diagonal in Figure~\ref{fig:figC}. This is expected because Equation~\eqref{eq:speedup-law} deliberately isolates the effect of \mymethod\ on the root cut loop and assumes that branching time and the remaining overhead are unchanged. In practice, the warm start may also affect these parts of the solve, for example through the primal start or through changes in the subsequent search. These deviations are therefore likely due to changes in parts of the solve that the model does not capture.

Overall, the trade-off in Equation~\eqref{eq:speedup-law} explains the main variation in \mymethod's performance: speedups occur when the warm start removes enough root cut loop time to outweigh the cost of constructing the surrogate.

\paragraph{Surrogate choice}
Finally, we address our third research question: how does the choice of surrogate affect \mymethod's performance? The results above suggest that two aspects matter: the surrogate should be inexpensive to solve, but it should also retain the information needed to generate useful cuts for the full problem. For \sstpp, TULIP-F illustrates the first point: solving only the first-stage block is sufficient to generate useful connectivity cuts, so keeping the stochastic coupling of TULIP-R mainly adds overhead. For SSTP, the second point is more important: the useful connectivity cuts are scenario-specific, making TULIP-S better suited to the relevant cut family than TULIP-R. The remaining case studies show that such a useful decomposition is not always available. All in all, these results suggest that surrogate choice involves a trade-off between computational cost and the ability to generate useful cuts for the full problem. The most effective surrogate is therefore problem-dependent.

\FloatBarrier

\section{Conclusion} \label{sec:conclusion}
We introduced \mymethod, a modular three-step warm-start framework for two-stage stochastic (mixed-)integer programs with an exponential number of cuts. Across four case studies, its performance is highly problem-dependent. \mymethod\ generally solves a similar number of instances to optimality as regular branch-and-cut, with some variants solving slightly more instances. On instances solved by all methods, the clearest runtime gains occur for the stochastic Steiner tree problems, where \mymethod\ accelerates regular branch-and-cut by a factor up to $2.56$, while providing little or no runtime benefit for windfarm, SCVRP, and SLR. On instances solved by no method, some \mymethod\ variants nevertheless attain smaller final optimality gaps, indicating greater progress toward optimality.

\mymethod\ accelerates branch-and-cut when the root cut loop time saved by the warm start exceeds the cost of constructing the surrogate, that is, when $(1-\alpha)f>\tau$. In our experiments, this occurs only for the stochastic Steiner tree problems, and most clearly for \sstpp. The \mymethod\ speedup model of Equation~\eqref{eq:speedup-law} explains the main differences in observed performance: a large baseline root cut loop fraction $f$, a large reduction in root cut loop time $1-\alpha$, and a small surrogate-overhead fraction $\tau$ favor a speedup. The choice of surrogate affects performance through both its computational cost and its ability to generate cuts that reduce the root cut loop work of the full problem.

These conclusions are based on four case studies, and their generality therefore remains to be established. Moreover, our experimental design favors instances that can be solved to optimality within the one-hour time limit, resulting in a relatively large share of easier instances. In such instances, the surrogate overhead accounts for a larger fraction of the total solution time and is therefore harder to compensate for. Further research should examine whether \mymethod\ becomes more beneficial on more challenging instances, for which the surrogate overhead may use a smaller fraction of the total solution time. Still, the strong results for the stochastic Steiner tree problem motivate further study of \mymethod\ for stochastic network design problems with dynamically separated connectivity constraints. A natural next step is to test whether the same speedup model continues to explain performance on other problem classes. A practical prospective procedure could first use a branch-and-cut run to measure $f$ and apply the baseline-only screening bound, and then, for promising instances, estimate $\tau$ from a surrogate-only pilot and $\alpha$ from a root-only warm-start pilot. Finally, it would be interesting to study whether the cut-reuse principle underlying \mymethod\ extends to feasibility cuts in integer L-shaped methods.

\section*{Acknowledgments}
We thank Ruurd Buijs and Markus Leitner for their insightful input, feedback, and fruitful discussions.

We thank SURF (\url{www.surf.nl}) for the support for using the National Supercomputer Snellius. This publication is part of the projects COCONUT and READINESS with project numbers 18938 and TWM.BL.019.002, respectively, which are financed by the Dutch Research Council (NWO).

\clearpage
\appendix
\raggedbottom %
\makeatletter
\setlength{\@fptop}{0pt}
\setlength{\@fpsep}{1.5\baselineskip}
\setlength{\@fpbot}{0pt plus 1fil}
\makeatother
\AddToHook{cmd/section/before}{\FloatBarrier}
\AddToHook{cmd/subsection/before}{\FloatBarrier}

\section{Stability results for TULIP-R} \label{sec:stability_tulip_r}
To assess the sensitivity of TULIP-R to the scenario-reduction seed, we consider its three runs $r \in \mathcal{R}_{\mathrm{seed}}$ separately. For $k \in \{0,1,2,3\}$, let
\[
    \mathcal{J}_{R}^{(k)} := \left\{ j \in \mathcal{J} : \text{TULIP-R solves instance $j$ to optimality for exactly $k$ runs in } \mathcal{R}_{\mathrm{seed}} \right\}.
\]
Thus,
\[
    \mathcal{J}_{R}^{(3)} = \mathcal{J}_{\mathrm{R}}^\star.
\]
Table~\ref{tab:tulipr-seed-robustness} reports the corresponding numbers of instances, aggregated over all scenario counts of each case study.

For instances in $\mathcal{J}_{R}^{(3)}$, we report the mean ratio between the longest and shortest seed runtimes in column ``Mean $T^{\max}/T^{\min}$'',
\[
    \frac{1}{|\mathcal{J}_{R}^{(3)}|} \sum_{j \in \mathcal{J}_{R}^{(3)}} \frac{ \max_{r \in \mathcal{R}_{\mathrm{seed}}} T^{(r)}_{j,\mathrm{TULIP\text{-}R}}}{\min_{r \in \mathcal{R}_{\mathrm{seed}}} T^{(r)}_{j,\mathrm{TULIP\text{-}R}}
    }.
\]
For instances in $\mathcal{J}_{R}^{(0)}$, we report the mean range between the largest and smallest final gaps across seeds in column ``Mean gap range (pp)'',
\[
    \frac{1}{|\mathcal{J}_{R}^{(0)}|} \sum_{j \in \mathcal{J}_{R}^{(0)}} \left( \max_{r \in \mathcal{R}_{\mathrm{seed}}} g^{(r)}_{j,\mathrm{TULIP\text{-}R}} - \min_{r \in \mathcal{R}_{\mathrm{seed}}} g^{(r)}_{j,\mathrm{TULIP\text{-}R}} \right),
\]
measured in percentage points (pp). As in Section~\ref{sec:performance_measures}, a run without an incumbent is assigned a final gap of $100\%$.

\begin{table}[!htbp]
\centering
\caption{Robustness of TULIP-R across the three scenario-reduction seeds.}
\label{tab:tulipr-seed-robustness}
\small
\begin{tabular}{lrrrrrrr}
\toprule
Problem
& $|\mathcal{J}|$
& $|\mathcal{J}_{R}^{(3)}|$
& $|\mathcal{J}_{R}^{(2)}|$
& $|\mathcal{J}_{R}^{(1)}|$
& $|\mathcal{J}_{R}^{(0)}|$
& Mean $T^{\max}/T^{\min}$
& Mean gap range (pp) \\
\midrule
\sstpp   & 77 & 74 & 1 & 0 & 2  & 1.22 & 0.0 \\
SSTP     & 77 & 77 & 0 & 0 & 0  & 1.05 & -- \\
SLR      & 80 & 62 & 4 & 7 & 7  & 1.25 & 16.5 \\
\slrp    & 60 & 48 & 2 & 0 & 10 & 1.38 & 28.4 \\
Windfarm & 18 & 13 & 0 & 0 & 5  & 1.04 & 0.0 \\
SCVRP    & 27 & 10 & 1 & 0 & 16 & 1.08 & 7.7 \\
\bottomrule
\end{tabular}
\end{table}
TULIP-R's ability to solve an instance to optimality is largely reproducible across seeds: all three seeds yield the same solve status for $96\%$ of the instances, increasing to $99\%$ when SLR and \slrp{} are excluded. Runtime variability is also limited: the mean ratio of the longest to the shortest seed runtimes ranges from $1.04$ to $1.38$ across the case studies. This variation is nevertheless comparable in magnitude to some of the reported speedups, so small differences in TULIP-R runtime should be interpreted with caution. SLR and \slrp{} are the most seed-sensitive case studies, particularly with respect to their final optimality gaps.

\section{Additional insight into impact warm-start cuts} \label{sec:figures_appendix}

To assess whether the transferred cuts strengthen the problem already at
the root node, we compare the first root-relaxation bounds reported by
Gurobi. Intuitively, we measure how much closer the first reported root
bound is to the final optimal objective under \mymethod\ than under regular
branch-and-cut. This can be seen as the change in the root-bound gap
to the optimum, rather than as a solver optimality gap, since no incumbent
for the full problem need be available at this point. For instance $j$ and method $m$, let $z^{\mathrm{log}}_{j,m}$ denote the objective reported in the first \textit{Root relaxation: objective} line of the final full-problem solve. Let $z^\star_{j,\mathrm{B\&C}}$ denote the final incumbent objective obtained by regular branch-and-cut, whose optimality is proven within the prescribed $0.01\%$ relative tolerance. We define the logged root-bound change as
\[
    \Delta^{\mathrm{log}}_{j,m} := 100\, \frac{ z^{\mathrm{log}}_{j,m} - z^{\mathrm{log}}_{j,\mathrm{B\&C}}}{|z^\star_{j,\mathrm{B\&C}}|}.
\]
Positive values indicate that method $m$ reduces the gap between the first reported root bound and the final optimal objective relative to regular branch-and-cut, whereas negative values indicate that this gap is larger. These logged objectives reflect the solver's root relaxation after MIP preprocessing and need not coincide with the optima of the static LP relaxations. The measure therefore captures a solver-level bound change rather than isolating the strengthening of a common LP relaxation by the injected cuts. Negative values are retained.

\begin{figure}[!htb]
    \centering
    \begin{subfigure}[t]{0.49\linewidth}
        \centering
        \includegraphics[width=\linewidth,trim=24bp 20bp 47bp 12bp,clip]{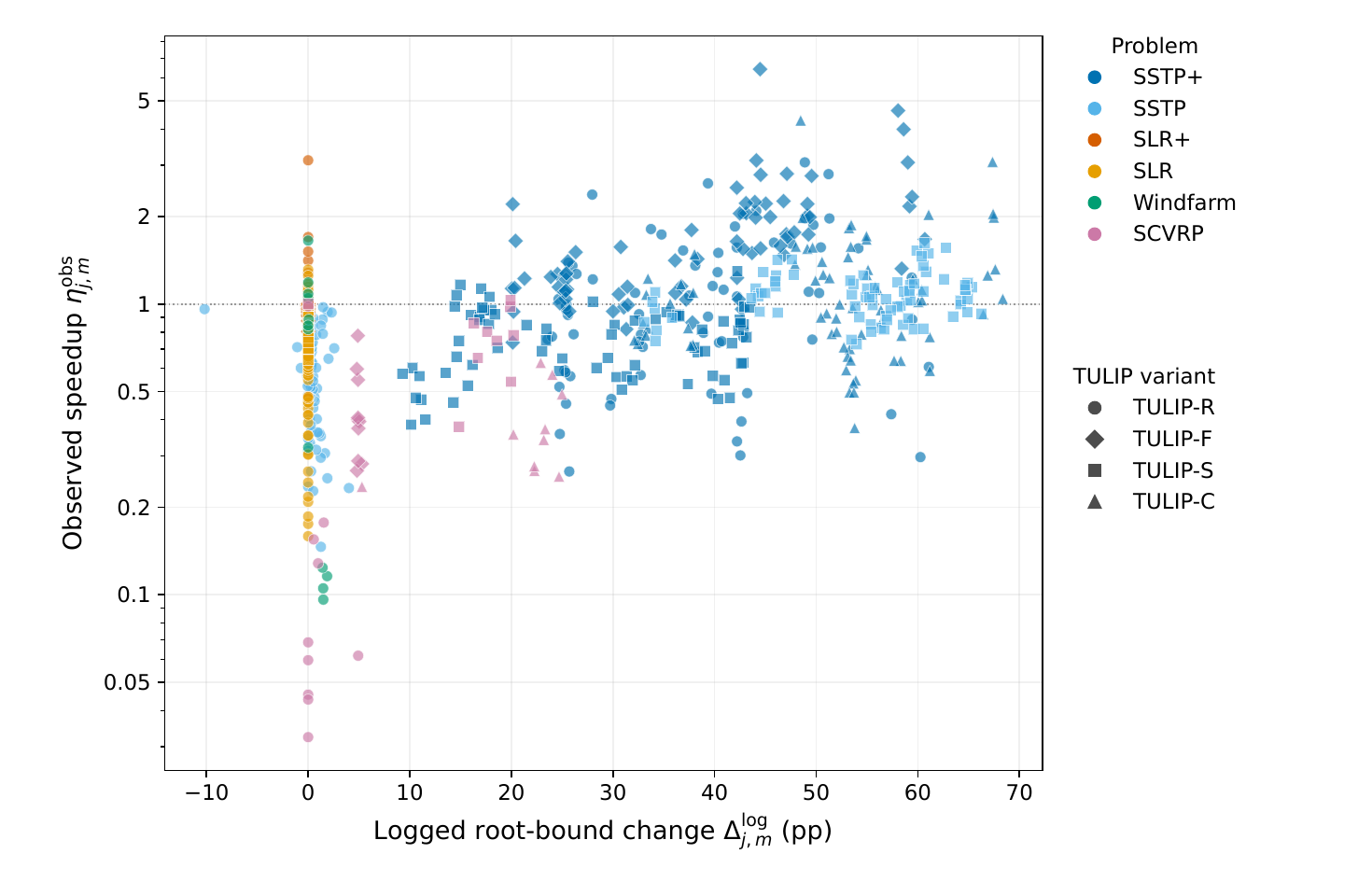}
        \subcaption{Observed speedup $\eta^{\mathrm{obs}}_{j,m}$ versus logged root-bound change $\Delta^{\mathrm{log}}_{j,m}$.}
        \label{fig:figA}
    \end{subfigure}\hfill
    \begin{subfigure}[t]{0.49\linewidth}
        \centering
        \includegraphics[width=\linewidth,trim=31bp 23bp 47bp 12bp,clip]{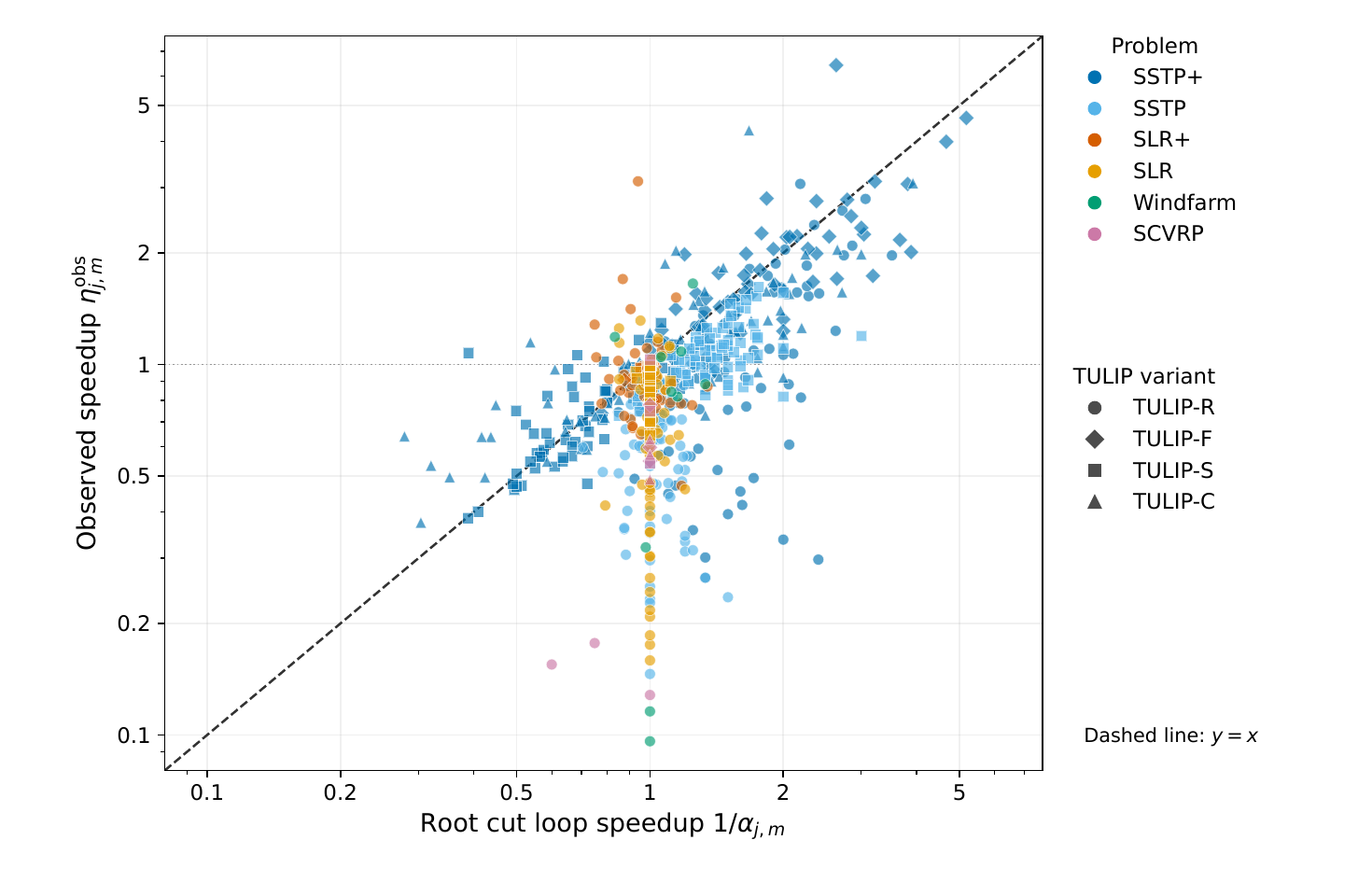}
        \subcaption{Observed speedup $\eta^{\mathrm{obs}}_{j,m}$ versus root cut loop speedup $1/\alpha_{j,m}$ across all case studies and applicable TULIP variants. The dashed line indicates $y=x$.}
        \label{fig:figB}
    \end{subfigure}

    \caption{Supporting empirical relationships between root-node effects and the observed speedup of \mymethod. Figure~\ref{fig:figA} relates the logged root-bound change to observed speedup, while Figure~\ref{fig:figB} relates the reduction in root cut loop time directly to the overall speedup.}
    \label{fig:root-bottleneck-accounting}
\end{figure}

Figure~\ref{fig:root-bottleneck-accounting} provides two diagnostics that support the mechanism underlying the \mymethod\ speedup model. While these relationships are not themselves the full model from Equation~\eqref{eq:speedup-law}, they help explain how the transferred cuts affect \mymethod's performance through the root node.

\FloatBarrier
Figure~\ref{fig:figA} shows that larger improvements in the first logged root-bound tend to coincide with larger observed speedups. This association is consistent with improved root-bounds contributing to faster solves. However, the logged values also reflect solver preprocessing, so the figure does not isolate the effect of the injected cuts on the static LP relaxation. The distinction between the stochastic Steiner tree problems and the remaining case studies is also visible: the largest root-bound improvements occur mostly for SSTP and \sstpp.

Figure~\ref{fig:figB} shows this mechanism more directly. Larger root cut loop speedups $1/\alpha_{j,m}$ are generally associated with larger overall speedups, whereas observations in which the warm start leaves the root cut loop unchanged or makes it slower tend to yield little or no overall benefit. Together, these figures provide supporting evidence for the interpretation used in Section~\ref{sec:why}: the usefulness of the transferred cuts is reflected primarily in their ability to reduce the computational effort at the root.

\section{Case study I: stochastic Steiner tree problem (SSTP)} \label{sec:case1}
\subsection{Problem}
We consider the SSTP, a two-stage stochastic version of the classical Steiner tree problem on a graph, following~\cite{leitner_decomposition_2018}. The problem represents network-design settings in which part of a network is installed before future connection requirements are known. In the standard SSTP, edges can be installed in the first stage at relatively low cost. These edges do not need to form a connected network. Uncertainty is represented by a finite set of scenarios, each with its own set of terminals. After a scenario is realized, additional edges can be purchased in the second stage at a higher cost to connect its terminals. The objective is to minimize the first-stage installation cost plus the expected second-stage purchase cost. Both stages use binary edge variables. Connectivity is enforced by an exponentially large set of connectivity cuts that are dynamically separated during branch-and-cut via a lazy constraint callback.

We also consider a variant in which the first-stage terminals must be connected by first-stage edges through Constraint~\eqref{eq:w1-cuts}. This variant is motivated by future-proof infrastructure design, such as ship pipe routing, where the network installed today must already connect the current system while allowing for future extensions and adaptations~\citep{markhorst2025future}. We use SSTP for the standard problem, in which first-stage connectivity is not required, and \sstpp\ for the extended version. The formulation in Appendix~\ref{sec:model1} describes \sstpp, and removing the first-stage connectivity cuts in~\eqref{eq:w1-cuts} gives the standard SSTP.

\subsection{Example}
To illustrate the SSTP, we use a small example based on~\cite[Figure 2]{markhorst_sailing_2024}. Consider the graph in Figure~\ref{fig:sstp_instance}, which shows the first-stage edge costs. The second-stage purchase costs are twice as high as the corresponding first-stage costs. The first-stage terminals are vertices $A$ and $D$, which must be connected using first-stage edges.

There are two scenarios, so $\mathcal{S}=\{1,2\}$. In scenario $1$, which occurs with probability $p^{(1)}=1-p^{(2)}$, no additional terminal appears. In scenario $2$, which occurs with probability $p^{(2)}$, vertex $B$ becomes an additional terminal that must be connected in the second stage.

Figure~\ref{fig:toy_example_sstp} shows the difference between the here-and-now-only and stochastic solutions. The here-and-now-only solution ignores the future scenarios and connects $A$ and $D$ using the cheapest edge $A$--$D$. The stochastic solution instead installs the first-stage edges $A$--$C$ and $C$--$D$. If scenario $2$ occurs, terminal $B$ can then be connected by purchasing only the second-stage edge $B$--$C$. For $p^{(2)}=0.4$, the stochastic solution has an expected cost of $2+0.4\cdot2=2.8$, compared with $1.5+0.4\cdot4=3.1$ for the here-and-now-only solution.

\begin{figure}[!htbp]
    \centering
    \begin{minipage}{0.45\linewidth}
        \centering
        \begin{tikzpicture}[node distance={15mm}, thick, square/.style={draw, regular polygon, regular polygon sides=4}]
            \node[circle, draw] (A) {A};
            \node[circle, draw] (C) [below of=A] {C};
            \node[circle, draw] (B) [left of=C] {B};
            \node[circle, draw] (D) [right of=C] {D};
            \draw (A) -- (C) node [midway, left] {1};
            \draw (A) -- (D) node [midway, above right] {1.5};
            \draw (B) -- (C) node [midway, below] {1};
            \draw (C) -- (D) node [midway, below] {1};
        \end{tikzpicture}
        \subcaption{Graph with the first-stage edge costs.}
        \label{fig:sstp_instance}
    \end{minipage}
    \vfill
    \begin{minipage}{0.45\linewidth}
        \centering
        \begin{tikzpicture}[node distance={15mm}, thick, square/.style={draw, regular polygon, regular polygon sides=4}]
            \node[circle, draw] (A) {A};
            \node[circle, draw] (C) [below of=A] {C};
            \node[circle, draw] (B) [left of=C] {B};
            \node[circle, draw] (D) [right of=C] {D};
            \draw (A) -- (C);
            \draw[blue, line width=3pt] (A) -- (D);
            \draw (B) -- (C);
            \draw (C) -- (D);
        \end{tikzpicture}
        \subcaption{Here-and-now-only solution.}
        \label{fig:sstp_do}
    \end{minipage}
    \hfill
    \begin{minipage}{0.45\linewidth}
        \centering
        \begin{tikzpicture}[node distance={15mm}, thick, square/.style={draw, regular polygon, regular polygon sides=4}]
            \node[circle, draw] (A) {A};
            \node[circle, draw] (C) [below of=A] {C};
            \node[circle, draw] (B) [left of=C] {B};
            \node[circle, draw] (D) [right of=C] {D};
            \draw[blue, line width=3pt] (A) -- (C);
            \draw (A) -- (D);
            \draw[blue, line width=3pt, dotted] (B) -- (C);
            \draw[blue, line width=3pt] (C) -- (D);
        \end{tikzpicture}
        \subcaption{Stochastic solution for $p^{(2)}=0.4$.}
        \label{fig:sstp_sp}
    \end{minipage}
    \caption{Small example illustrating the SSTP, based on~\cite[Figure 2]{markhorst_sailing_2024}. Highlighted edges are installed: thick solid lines denote first-stage edges and thick dotted lines denote second-stage (purchased) edges; thin lines are the underlying graph.}
    \label{fig:toy_example_sstp}
\end{figure}

\subsection{Model} \label{sec:model1}
The formulation follows the modeling setup of our earlier work on future-proof ship pipe routing~\citep{markhorst2025future}, with the main difference that we consider stochastic Steiner \textit{trees} rather than Steiner \textit{forests}.

\paragraph{Sets and parameters}
Let $G = (V, E)$ be an undirected graph and let $A$ be the set of directed arcs, where each undirected edge $\{i,j\} \in E$ gives rise to arcs $(i,j)$ and $(j,i)$. Let $T \subseteq V$ be the set of first-stage terminals, $\mathcal{S}$ the set of scenarios with probabilities $p^{(s)} > 0$, $c_e$ the first-stage installation cost of edge $e$, $q_e^{(s)}$ the second-stage purchase cost of edge $e$ in scenario $s$, and $\tau_v^{(s)} \in \{0,1\}$ an indicator equal to $1$ if node $v$ is a terminal in scenario $s$. The relation between the two terminal levels is a modeling choice of our study: The first-stage terminal set is taken to be the terminal set of the \emph{first} scenario, $T = \{v : \tau_v^{(1)} = 1\}$, and each scenario $s$ requires connectivity only among its own terminals $\{v : \tau_v^{(s)} = 1\}$, which need not contain $T$. In the connectivity cuts, $r$ denotes an arbitrary fixed first-stage terminal in $T$ serving as root, and $r^{(s)}$ an arbitrary fixed terminal of scenario $s$; our implementation uses the smallest-index terminal in both cases. For each set $U \subseteq V$, let $\delta^-(U) = \{(u,v)\in A : u \notin U,\ v \in U\}$ and $\delta^+(U) = \{(u,v)\in A : u \in U,\ v \notin U\}$ denote the sets of arcs entering and leaving $U$, respectively; for a singleton set $U = \{v\}$, we simply write $\delta^-(v)$ and $\delta^+(v)$.

\paragraph{Decision variables}
To model connectivity, we complement the undirected edge-installation decisions with directed routing variables. These variables direct selected installed edges into a routing structure from the root to the terminals. The model uses four families of binary variables. The first-stage decisions are the edge-install variables $x_e \in \{0,1\}$ for $e \in E$, equal to $1$ when edge $e$ is installed in the first stage, together with the routing variables $w^1_a \in \{0,1\}$ for $a \in A$, equal to $1$ when arc $a$ is used in the first-stage routing. The second-stage decisions are taken per scenario: The purchase variables $y_e^{(s)} \in \{0,1\}$ for $e \in E$ and $s \in \mathcal{S}$, equal to $1$ when edge $e$ is additionally bought in scenario $s$, and the routing variables $w^{2,(s)}_a \in \{0,1\}$ for $a \in A$ and $s \in \mathcal{S}$, equal to $1$ when arc $a$ is used in the routing of scenario $s$.

\paragraph{Formulation}
\begin{subequations}
\begin{alignat}{3}
  \min \quad
    & \sum_{e \in E} c_e\, x_e
    + \sum_{s \in \mathcal{S}} p^{(s)} \sum_{e \in E} q_e^{(s)}\, y_e^{(s)}
    \label{eq:obj} \\[6pt]
  \text{s.t.} \quad
    & w^1_{(i,j)} + w^1_{(j,i)} \le x_{\{i,j\}},
    && \quad \{i,j\} \in E,
    \label{eq:w1-install} \\
    & \sum_{a \in \delta^-(v)} w^1_a \le 1,
    && \quad v \in V,
    \label{eq:w1-indegree} \\
    & \sum_{a \in \delta^-(v)} w^1_a \le \sum_{a \in \delta^+(v)} w^1_a,
    && \quad v \in V \setminus T,
    \label{eq:w1-steiner} \\
    & \sum_{a \in \delta^+(U)} w^1_a \ge 1,
    && \quad U \subset V,\; r \in U,\; t \in V \setminus U,\; t \in T,
    \label{eq:w1-cuts} \\[6pt]
    & w^{2,(s)}_{(i,j)} + w^{2,(s)}_{(j,i)} \le x_{\{i,j\}} + y^{(s)}_{\{i,j\}},
    && \quad \{i,j\} \in E,\; s \in \mathcal{S},
    \label{eq:w2-install} \\
    & w^{2,(s)}_{(i,j)} + w^{2,(s)}_{(j,i)} \le 1,
    && \quad \{i,j\} \in E,\; s \in \mathcal{S},
    \label{eq:w2-direction} \\
    & \sum_{a \in \delta^-(v)} w^{2,(s)}_a \le 1,
    && \quad v \in V,\; s \in \mathcal{S},
    \label{eq:w2-indegree} \\
    & (1 - \tau_v^{(s)})\sum_{a \in \delta^-(v)} w^{2,(s)}_a
      \le (1 - \tau_v^{(s)})\sum_{a \in \delta^+(v)} w^{2,(s)}_a,
    && \quad v \in V,\; s \in \mathcal{S},
    \label{eq:w2-steiner} \\
    & \sum_{a \in \delta^+(U)} w^{2,(s)}_a \ge 1,
    && \quad \begin{aligned}[t]
         & U \subset V,\; r^{(s)} \in U,\\
         & t \in V \setminus U,\; \tau_t^{(s)} = 1,\; s \in \mathcal{S},
       \end{aligned}
    \label{eq:w2-cuts} \\[6pt]
    & x_e,\; y_e^{(s)} \in \{0,1\},
    && \quad e \in E,\; s \in \mathcal{S}, \\
    & w^1_a,\; w^{2,(s)}_a \in \{0,1\},
    && \quad a \in A,\; s \in \mathcal{S}.
\end{alignat}
\end{subequations}

Constraints~\eqref{eq:w1-install} and~\eqref{eq:w2-install} ensure that arcs can be routed only along installed edges, either in the first stage or as additional purchases in the second stage. Constraints~\eqref{eq:w1-indegree} and~\eqref{eq:w2-indegree} impose that each node receives at most one incoming arc, enforcing a tree structure on the routing. Constraints~\eqref{eq:w1-steiner} and~\eqref{eq:w2-steiner} encode the Steiner-node pass-through condition: if the routing enters a non-terminal node, it must also leave that node, so the routing cannot terminate at a Steiner node; in the second-stage constraint, the indicator $(1 - \tau_v^{(s)})$ activates this condition only for nodes that are not terminals in scenario $s$. Constraint~\eqref{eq:w2-direction} additionally prevents both directions of the same edge from being active simultaneously in any scenario. Finally, constraints~\eqref{eq:w1-cuts} and~\eqref{eq:w2-cuts} are connectivity cuts, added lazily via a callback: every subset $U \subset V$ containing the root but not some terminal must have at least one outgoing arc active in the routing, guaranteeing that all terminals are reachable in both the first stage and in every scenario. Constraints~\eqref{eq:w1-cuts} and~\eqref{eq:w2-cuts} are connectivity-cut inequalities, which are added lazily via a callback. To separate these inequalities, we use graph-theoretic minimum-cut problems. For each terminal $t$, the current routing-variable values are used as arc capacities and a minimum $(r,t)$-cut is computed. If its capacity is smaller than one, the corresponding connectivity-cut inequality is violated and added to the model. Repeating this separation procedure at the following callback solutions may generate different connectivity cuts. For the backward inequalities, which are not required to enforce connectivity but are added to strengthen the LP relaxation and improve computational performance, we reverse the arc capacities and similarly solve a minimum $(t,r)$-cut problem. Any violated cut of either type is added as a lazy constraint during branch-and-cut, so only cuts that are actually violated at the current node are imposed.

\subsection{Experiments} \label{sec:experiments1}
We use the 11 graphs K100 and K100.1--K100.10 from the K100 family of the SteinLib library~\citep{koch2001steinlib}. Their main characteristics are summarized in Table~\ref{tab:steiner-instances}. The first-stage edge costs are taken from SteinLib and are the same across scenarios. Each scenario has its own terminal set and second-stage edge costs.

For each graph, we consider scenario counts $|\mathcal{S}|\in\{50,75,100,150,200,300,500\}$. For \sstpp{}, we compare regular branch-and-cut with TULIP-R, TULIP-F, TULIP-S, and TULIP-C. For SSTP, only TULIP-R and TULIP-S are considered, since there are no first-stage connectivity cuts. The second-stage surrogates assume no preinstalled edges ($x_e=0$ for all $e\in E$), and their construction and solution times are included in $t_{\mathrm{sur}}$. The remaining experimental settings are described in Section~\ref{sec:experimental_design}.

\begin{table}[!htbp]
\caption{Basic properties of the $11$ graphs used in the SSTP study (K100 and K100.1--K100.10): number of vertices $|V|$, edges $|E|$, and first-stage terminals $|K|$ (computed from the committed instance files; the graph and the first-stage terminal set are shared across all scenario counts of an instance).}
\label{tab:steiner-instances}
\centering
\begin{tabular}{l r rrr rrr rrr}
\toprule
Dataset & Inst [\#] & \multicolumn{3}{l}{$|V|$} & \multicolumn{3}{l}{$|E|$} & \multicolumn{3}{l}{$|K|$} \\
\cmidrule(r){3-5} \cmidrule(lr){6-8} \cmidrule(l){9-11}
 & & Min & Avg & Max & Min & Avg & Max & Min & Avg & Max \\
\midrule
K100     & 11 & 22   & 30.5   & 45   & 64   & 114.9  & 191    & 3 & 5.5 & 7 \\
\bottomrule
\end{tabular}
\end{table}

\subsection{Results} \label{sec:case1_results}

\begin{table}[!htbp]
\centering
\caption{Number of SSTP instances solved and mean speedup split by scenario count. Bold values indicate the best performance within each scenario count, and ties are all highlighted.}
\label{tab:solved-speedup-steiner-nofs}
\small
\begin{tabular}{rrrcccccc}
\toprule
&&& \multicolumn{2}{c}{B\&C} & \multicolumn{2}{c}{TULIP-R} & \multicolumn{2}{c}{TULIP-S} \\
\cmidrule(lr){4-5}\cmidrule(lr){6-7}\cmidrule(lr){8-9}
$|\mathcal{S}|$
& $|\mathcal{J}|$
& $|\mathcal{J}^{\cap}|$
& $|\mathcal{J}_{\mathrm{B\&C}}^\star|$
& $\bar{T}_{\mathrm{B\&C}}$ (s)
& $|\mathcal{J}_{R}^\star|$
& $\bar{\eta}_{R}$
& $|\mathcal{J}_{S}^\star|$
& $\bar{\eta}_{S}$ \\
\midrule
$50$  & 11 & 11 & \textbf{11\,(100\%)} & \textbf{9.7} & \textbf{11\,(100\%)} & 0.37 & \textbf{11\,(100\%)} & 0.99 \\
$75$  & 11 & 11 & \textbf{11\,(100\%)} & \textbf{17}  & \textbf{11\,(100\%)} & 0.45 & \textbf{11\,(100\%)} & 0.99 \\
$100$ & 11 & 11 & \textbf{11\,(100\%)} & 24  & \textbf{11\,(100\%)} & 0.51 & \textbf{11\,(100\%)} & \textbf{1.06} \\
$150$ & 11 & 11 & \textbf{11\,(100\%)} & 44  & \textbf{11\,(100\%)} & 0.62 & \textbf{11\,(100\%)} & \textbf{1.05} \\
$200$ & 11 & 11 & \textbf{11\,(100\%)} & 78  & \textbf{11\,(100\%)} & 0.68 & \textbf{11\,(100\%)} & \textbf{1.07} \\
$300$ & 11 & 11 & \textbf{11\,(100\%)} & 156 & \textbf{11\,(100\%)} & 0.74 & \textbf{11\,(100\%)} & \textbf{1.17} \\
$500$ & 11 & 11 & \textbf{11\,(100\%)} & 388 & \textbf{11\,(100\%)} & 0.83 & \textbf{11\,(100\%)} & \textbf{1.18} \\
\bottomrule
\end{tabular}
\end{table}

\begin{figure}[!htbp]
    \centering
    \includegraphics[width=0.65\linewidth]{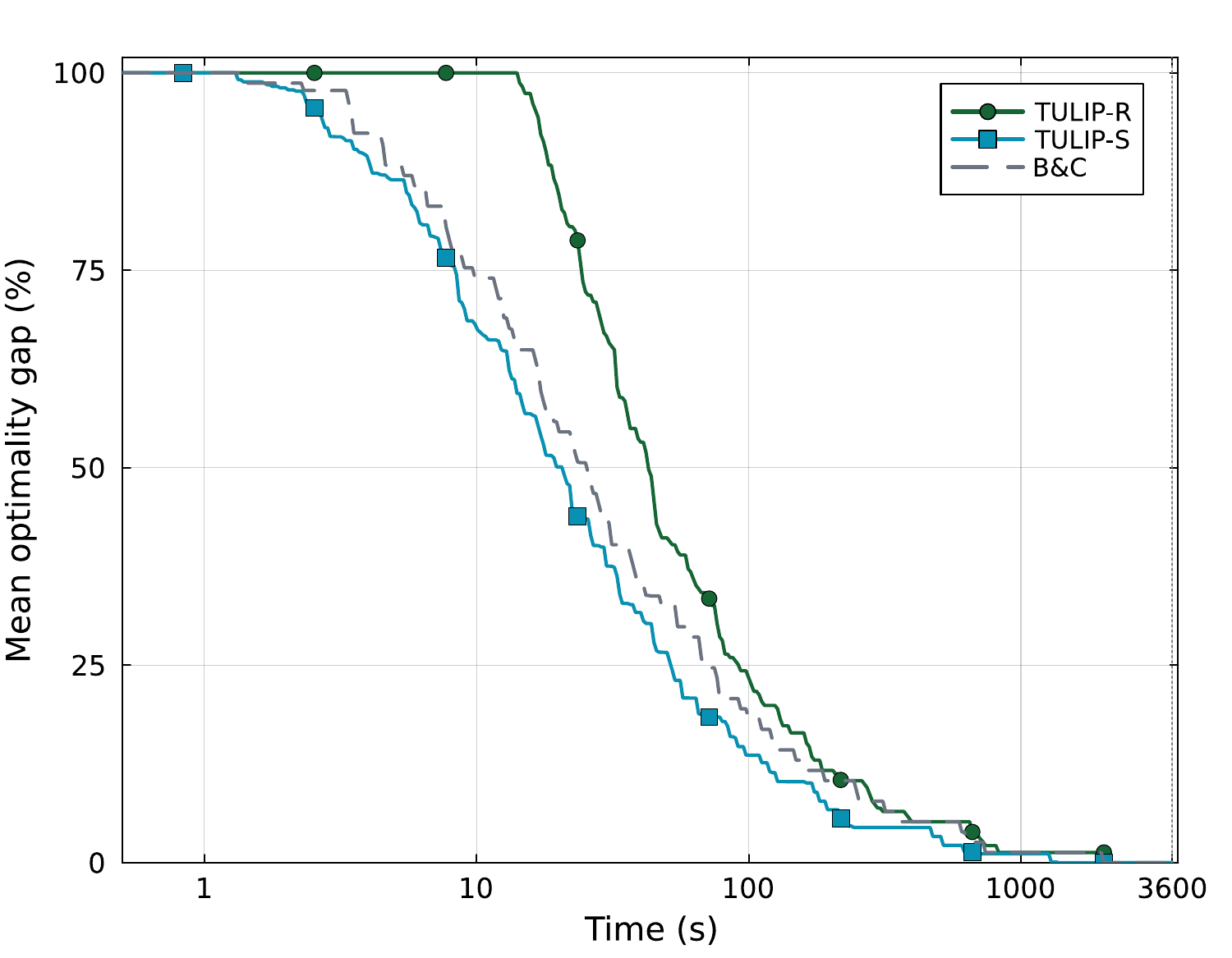}
    \caption{Mean optimality gap over time for SSTP.}
    \label{fig:steiner_nofs}
\end{figure}

\begin{figure}[!htb]
    \centering
    \begin{subfigure}[t]{0.49\linewidth}
        \centering
        \includegraphics[width=\linewidth,trim=43bp 34bp 173bp 2bp,clip]{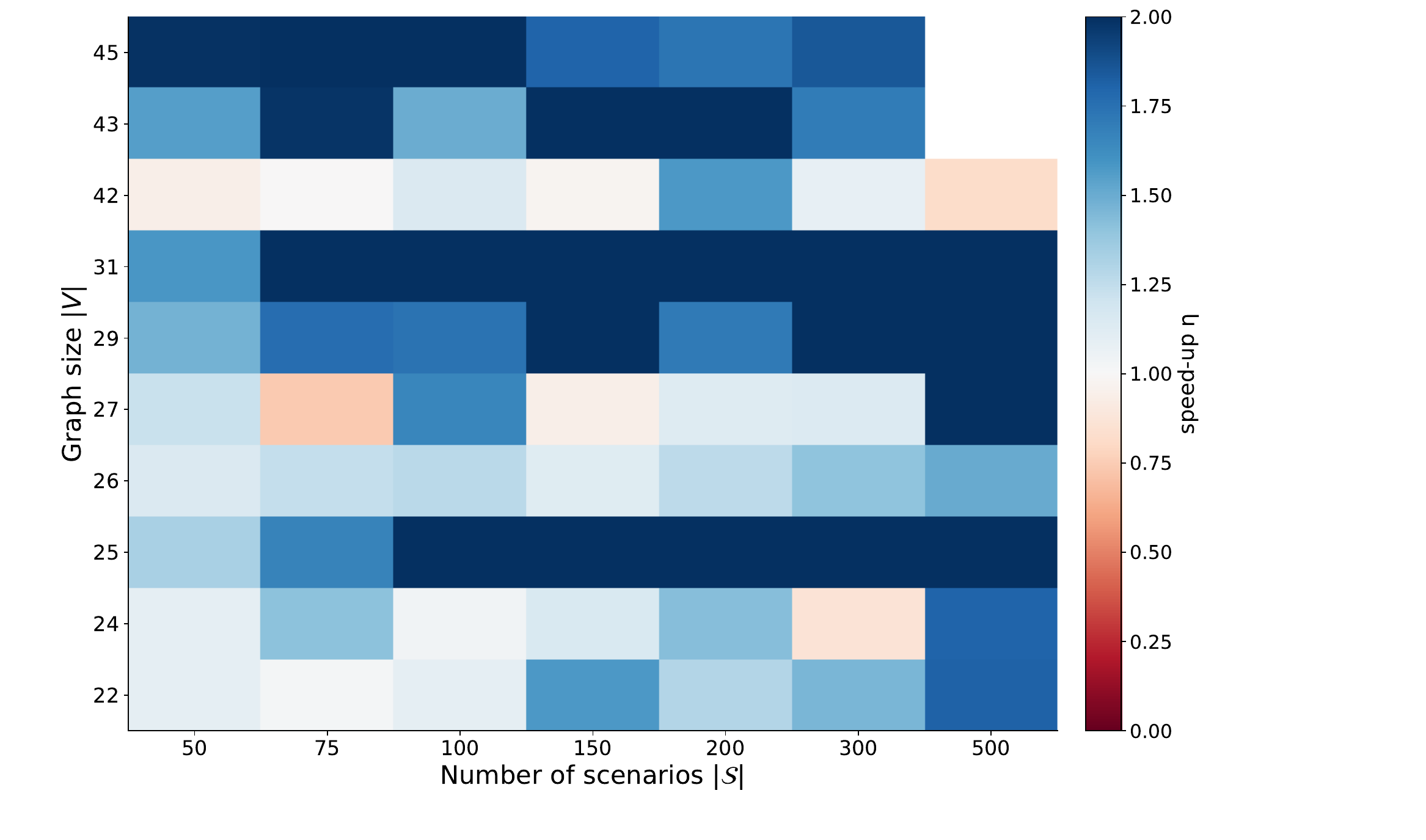}
        \subcaption{\sstpp\ with TULIP-F.}
        \label{fig:steiner-nodes-full}
    \end{subfigure}\hfill
    \begin{subfigure}[t]{0.49\linewidth}
        \centering
        \includegraphics[width=\linewidth,trim=43bp 34bp 173bp 2bp,clip]{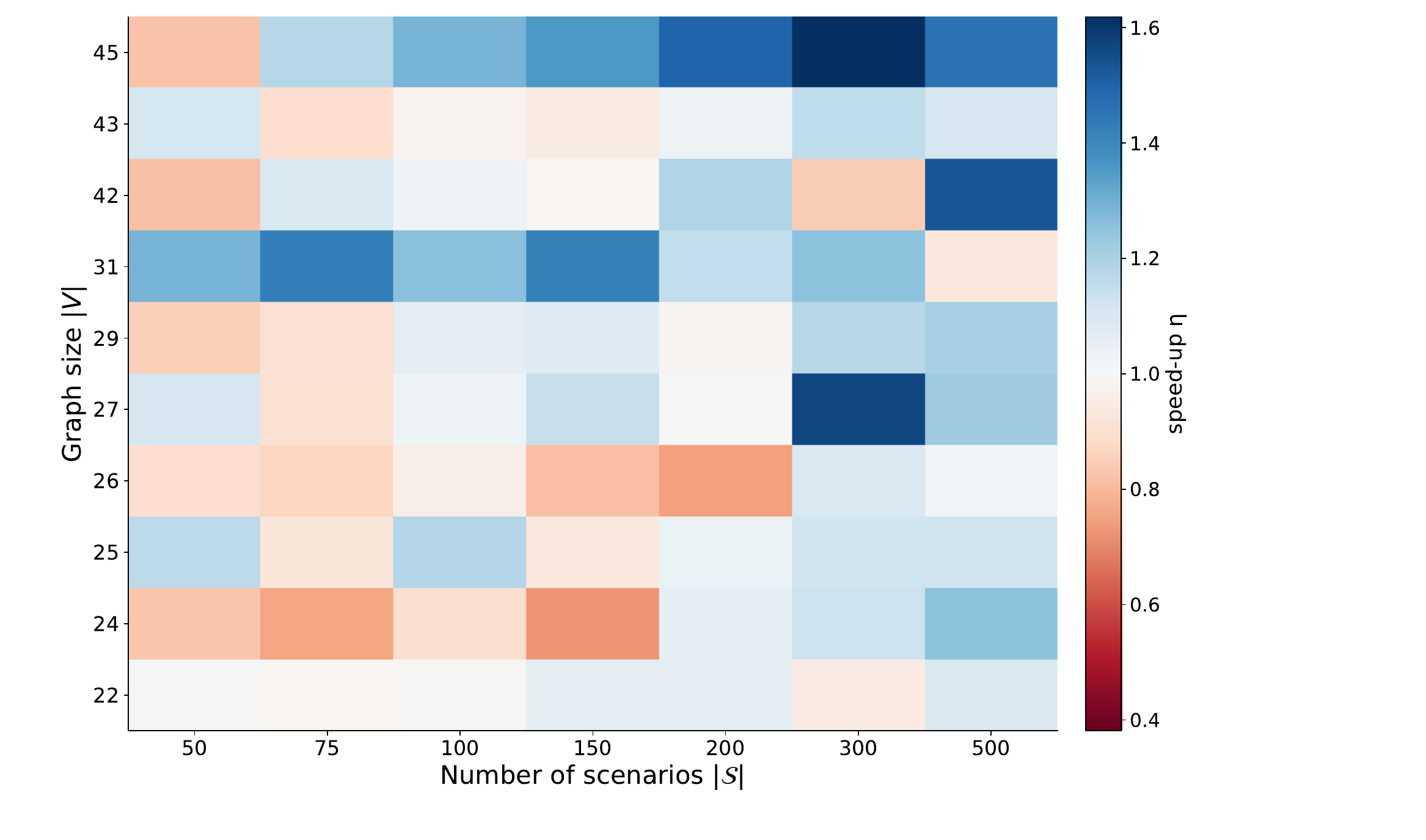}
        \subcaption{SSTP with TULIP-S.}
        \label{fig:steiner-nodes-nofs}
    \end{subfigure}
    \caption{Mean observed speedup over regular branch-and-cut by graph size $|V|$ and scenario count $|\mathcal{S}|$. Values above one indicate that TULIP is faster than regular branch-and-cut.}
    \label{fig:steiner-nodes}
\end{figure}

\FloatBarrier
\paragraph{Effect of graph size}
Figure~\ref{fig:steiner-nodes} reports the mean observed speedup by graph size and scenario count for \sstpp\ with TULIP-F and SSTP with TULIP-S.

\FloatBarrier
\paragraph{SAA benchmark}
To provide an approximate-solution benchmark for \mymethod, we draw a subset of scenarios of size $|\mathcal{S}'| = \lceil 0.1 | \mathcal{S} | \rceil$, solve the resulting SAA problem to optimality, and evaluate the obtained first-stage decision \emph{out-of-sample} on the full scenario set $\mathcal{S}$. We apply this procedure to both SSTP and \sstpp. Each experiment is repeated over three replications, and the out-of-sample objective is compared with the corresponding full-scenario optimum $z^\star$ obtained by branch-and-cut.

Unlike the methods in $\mathcal{M}$, SAA does not solve the full-scenario problem to optimality; instead, it returns an approximate first-stage solution based on a smaller sampled problem. We therefore evaluate SAA separately, reporting both its computational time and the out-of-sample quality of the resulting solution. This comparison provides a reference for the trade-off between the computational savings from solving a smaller problem and the loss in solution quality relative to the full-scenario optimum.

Let $\mathcal{R}_{\mathrm{SAA}}$ denote the set of SAA replications and let $T_{j,\mathrm{SAA}}^{(r)}$ denote the time required to solve replication $r$ for instance $j$. We define the mean SAA runtime for instance $j$ as
\[
    T_{j,\mathrm{SAA}} := \frac{1}{|\mathcal{R}_{\mathrm{SAA}}|} \sum_{r\in \mathcal{R}_{\mathrm{SAA}}} T_{j,\mathrm{SAA}}^{(r)}.
\]
Let $z_j^\star$ denote the full-scenario optimum and $z_{j,\mathrm{oos}}^{(r)}$ the objective obtained by evaluating the first-stage solution from replication $r$ on the full scenario set. The corresponding mean out-of-sample gap is
\[
    \delta_{j,\mathrm{SAA}} := \frac{ \frac{1}{|\mathcal{R}_{\mathrm{SAA}}|} \sum_{r\in \mathcal{R}_{\mathrm{SAA}}} z_{j,\mathrm{oos}}^{(r)} - z_j^\star }{|z_j^\star|} \cdot 100\%.
\]
Over a set of instances $\mathcal{J}$, we report the means
\[
    \bar T_{\mathrm{SAA}} := \frac{1}{|\mathcal{J}|} \sum_{j\in\mathcal{J}}T_{j,\mathrm{SAA}}, \qquad \bar\delta_{\mathrm{SAA}} := \frac{1}{|\mathcal{J}|} \sum_{j\in\mathcal{J}}\delta_{j,\mathrm{SAA}}.
\]
For the SAA gap, we additionally report the standard deviation across instances.

\begin{table}[!htbp]
\centering
\caption{Out-of-sample SAA results by scenario count for SSTP and \sstpp. Here, $|\mathcal{J}|$ is the number of instances with a certified regular branch-and-cut optimum used as the reference value.}
\label{tab:SAA-results-steiner}
\small
\begin{tabular}{rcccccc}
\toprule
& \multicolumn{3}{c}{SSTP} & \multicolumn{3}{c}{\sstpp} \\
\cmidrule(lr){2-4}\cmidrule(lr){5-7}
$|\mathcal{S}|$ & $|\mathcal{J}|$ & $\bar{\delta}_{\mathrm{SAA}}$ (\%) & $\bar{T}_{\mathrm{SAA}}$ (s) & $|\mathcal{J}|$ & $\bar{\delta}_{\mathrm{SAA}}$ (\%) & $\bar{T}_{\mathrm{SAA}}$ (s) \\
\midrule
$50$  & 11 & $2.06\pm1.54$ & $2.8$ & 11 & $1.09\pm0.65$ & $3.5$ \\
$75$  & 11 & $1.21\pm1.23$ & $3.2$ & 11 & $0.43\pm0.56$ & $4.5$ \\
$100$ & 11 & $1.28\pm1.14$ & $3.6$ & 11 & $0.62\pm0.57$ & $6.3$ \\
$150$ & 11 & $0.96\pm0.93$ & $4.3$ & 10 & $0.25\pm0.45$ & $5.6$ \\
$200$ & 11 & $0.38\pm0.49$ & $5.2$ & 11 & $0.13\pm0.18$ & $11.6$ \\
$300$ & 11 & $0.53\pm0.52$ & $7.9$ & 9 & $0.15\pm0.19$ & $6.8$ \\
$500$ & 11 & $0.06\pm0.11$ & $11.7$ & 9 & $0.03\pm0.06$ & $11.9$ \\
\bottomrule
\end{tabular}
\end{table}

Table~\ref{tab:SAA-results-steiner} reports the out-of-sample optimality gaps of the SAA first-stage decisions for both stochastic Steiner tree formulations. The gaps generally decrease as the scenario count grows. For \sstpp, the mean gap is already below $1\%$ at $|\mathcal{S}|\ge 75$; for SSTP, it falls from $2.06\%$ at $|\mathcal{S}|=50$ to $0.06\%$ at $|\mathcal{S}|=500$.

\FloatBarrier
\Needspace{24\baselineskip}
\paragraph{Expected value of perfect information}
Table~\ref{tab:evpi-steiner-full} reports the EVPI per K100 instance at $|\mathcal{S}|=50$. The mean EVPI is $11{,}796$, while the mean instance-level EVPI percentage is $6.8\%$, confirming that accounting for stochasticity in the case study is important. The value of the stochastic solution is omitted because it is degenerate under the binary-terminal uncertainty.

\begin{table}[H]
\centering
\caption{EVPI on the K100 \sstpp\ family, $|\mathcal{S}|=50$. $z^\star$ is the stochastic optimum and WS$=z^\star-\text{EVPI}$ the wait-and-see value. Values are rounded independently to the nearest integer, so the displayed identity may differ by one unit.}
\label{tab:evpi-steiner-full}
\begin{tabular}{lrrrr}
\toprule
Instance & $z^\star$ & WS & EVPI & EVPI (\%) \\
\midrule
K100 & $224{,}131$ & $209{,}548$ & $14{,}583$ & $6.5$ \\
K100.1 & $165{,}430$ & $152{,}254$ & $13{,}176$ & $8.0$ \\
K100.2 & $126{,}893$ & $115{,}026$ & $11{,}868$ & $9.4$ \\
K100.3 & $129{,}877$ & $115{,}965$ & $13{,}912$ & $10.7$ \\
K100.4 & $141{,}394$ & $131{,}759$ & $9{,}635$ & $6.8$ \\
K100.5 & $220{,}780$ & $211{,}334$ & $9{,}445$ & $4.3$ \\
K100.6 & $183{,}270$ & $174{,}102$ & $9{,}168$ & $5.0$ \\
K100.7 & $184{,}259$ & $179{,}233$ & $5{,}027$ & $2.7$ \\
K100.8 & $246{,}218$ & $230{,}925$ & $15{,}292$ & $6.2$ \\
K100.9 & $146{,}414$ & $136{,}971$ & $9{,}444$ & $6.5$ \\
K100.10 & $213{,}794$ & $195{,}590$ & $18{,}204$ & $8.5$ \\
\midrule
\textbf{Mean} & $180{,}224$ & $168{,}428$ & $11{,}796$ & $6.8$ \\
\bottomrule
\end{tabular}
\end{table}

\section{Case study II: stochastic windfarm problem} \label{sec:case2}

\subsection{Problem}
In this section, we develop a stochastic version of the integrated windfarm layout and cable routing problem from \cite{pedersen2026integrated}. In the stochastic version considered here, there is one substation $r$ and a set of candidate turbine locations $P$. In the first stage, we decide which turbines to build and how to connect the built turbines to the substation through a cable network, while respecting a minimum distance between selected turbine locations. Building and connecting turbines incurs a first-stage investment cost consisting of cable costs and location-dependent foundation costs. Wind conditions are uncertain and are represented by a finite scenario set $\mathcal{S}$. In scenario $s$, each built turbine produces a scenario-dependent amount of power, while wake interference between built turbines may reduce the total production. If the resulting production in the second stage falls below the contracted quota $Q$, a shortfall penalty is incurred. The objective is to minimize the first-stage investment cost plus the expected shortfall penalty. The first stage can therefore be viewed as a Steiner tree problem in which the turbine locations to be connected are themselves decisions.

\subsection{Example}
To illustrate the windfarm problem, consider the small instance in Figure~\ref{fig:wf_instance}. There is one substation, denoted by root $r$, and two candidate turbine locations, $1$ and $2$. For simplicity, we aggregate the turbine installation and cable costs. Building a turbine and connecting it to the substation through the cable network has a first-stage cost of $3$. The two locations are far enough apart that both turbines can be built. The contracted quota is $Q=10$, and each unit of shortfall has a penalty cost of $\rho=3$.

There are two equally likely wind scenarios, $\mathcal{S}=\{1,2\}$. In the windy scenario, each built turbine produces $16$ units of power. In the calm scenario, each produces $4$ units. If both turbines are built, wake interference reduces their combined production by $4$ units in the windy scenario but does not in the calm scenario.

The deterministic expected-value problem uses the average production of each turbine, which is $10$ units. It therefore builds one turbine, as shown in Figure~\ref{fig:wf_do}. The first-stage investment cost is $3$. Under the actual scenarios, this turbine meets the quota in the windy scenario but has a shortfall of $6$ units in the calm scenario. Its expected cost is therefore
\[
    3 + \tfrac{1}{2}\cdot 3\cdot 6 = 12.
\]

The stochastic optimization solution instead builds both turbines, as shown in Figure~\ref{fig:wf_sp}. The first-stage investment cost is then $6$. The two turbines meet the quota in the windy scenario and have a shortfall of only $2$ units in the calm scenario. The expected cost is therefore
\[
    6 + \tfrac{1}{2}\cdot 3\cdot 2 = 9.
\]
Hence, taking the possible wind scenarios into account makes it optimal to build both turbines.

\begin{figure}[!htbp]
    \centering
    \begin{minipage}{0.45\linewidth}
        \centering
        \begin{tikzpicture}[thick]
            \node[draw, rectangle, minimum size=6mm] (S) at (0,0) {$r$};
            \node[circle, draw] (T1) at (2,0.8) {$1$};
            \node[circle, draw] (T2) at (2,-0.8) {$2$};
            \draw (S) -- (T1) node [midway, above left] {3};
            \draw (S) -- (T2) node [midway, below left] {3};
        \end{tikzpicture}
        \subcaption{Substation $r$ and two candidate turbines with cable costs.}
        \label{fig:wf_instance}
    \end{minipage}
    \vfill
    \begin{minipage}{0.45\linewidth}
        \centering
        \begin{tikzpicture}[thick]
            \node[draw, rectangle, minimum size=6mm] (S) at (0,0) {$r$};
            \node[circle, draw] (T1) at (2,0.8) {$1$};
            \node[circle, draw] (T2) at (2,-0.8) {$2$};
            \draw[blue, line width=3pt] (S) -- (T1);
            \draw (S) -- (T2);
        \end{tikzpicture}
        \subcaption{Deterministic solution: build one turbine.}
        \label{fig:wf_do}
    \end{minipage}
    \hfill
    \begin{minipage}{0.45\linewidth}
        \centering
        \begin{tikzpicture}[thick]
            \node[draw, rectangle, minimum size=6mm] (S) at (0,0) {$r$};
            \node[circle, draw] (T1) at (2,0.8) {$1$};
            \node[circle, draw] (T2) at (2,-0.8) {$2$};
            \draw[blue, line width=3pt] (S) -- (T1);
            \draw[blue, line width=3pt] (S) -- (T2);
        \end{tikzpicture}
        \subcaption{Stochastic solution: build both turbines.}
        \label{fig:wf_sp}
    \end{minipage}
    \caption{Small example illustrating the windfarm problem. A built turbine is connected to the substation $r$ by an installed cable (thick blue).}
    \label{fig:toy_example_wf}
\end{figure}

\subsection{Model}
Figure~\ref{fig:windfarm-transformation} illustrates the modeling concept. The first-stage problem can be viewed as a Steiner tree problem on the transformed graph, in which the substation $r$ must be connected to all dummy terminals in $D$. The candidate turbine locations in $P$ can serve as Steiner nodes. For each candidate location $i\in P$, the dummy terminal $i+n\in D$ can be connected either through $(i,i+n)$, which represents building turbine $i$, or through $(r,i+n)$, which represents not building it. Exactly one of these two options is selected for every candidate location. The selected turbines determine the power production in each scenario, which is reduced by wake interference and compared with the required quota $Q$. Solid arcs in Figure~\ref{fig:windfarm-transformation} indicate selected arcs, whereas dashed arcs indicate available but unselected arcs.

\begin{figure}[t]
    \centering
    \begin{tikzpicture}[
        >=stealth,
        node distance=1.8cm and 2.2cm,
        every node/.style={font=\small},
        root/.style={circle, draw, minimum size=7mm},
        turbine/.style={circle, draw, minimum size=7mm},
        dummy/.style={rectangle, draw, minimum size=7mm},
        selected/.style={->, very thick},
        available/.style={->, thin, dashed},
        available2/.style={-, thin, dashed}
    ]

    \node[root] (r) {$r$};

    \node[turbine, right=2.6cm of r, yshift=1.2cm] (i) {$i$};
    \node[turbine, right=2.6cm of r, yshift=-1.2cm] (j) {$j$};

    \node[dummy, right=2.6cm of i] (ip) {$i+n$};
    \node[dummy, right=2.6cm of j] (jp) {$j+n$};

    \draw[selected] (r) -- (i);
    \draw[selected] (i) -- (ip);
    \draw[selected] (r) -- (jp);

    \draw[available] (r) -- (ip);
    \draw[available] (j) -- (jp);
    \draw[available] (r) -- (j);
    \draw[available2] (i) -- (j);

    \node[above=3mm of r] {substation};
    \node[above=3mm of i] {candidate location};
    \node[above=3mm of ip] {dummy terminal};

    \node[right=4mm of ip, align=left]
        {$i$ is built:\\[-1mm] $(i,i+n)$ selected};

    \node[right=4mm of jp, align=left]
        {$j$ is not built:\\[-1mm] $(r,j+n)$ selected};

    \end{tikzpicture}
    \caption{Illustration of the transformed windfarm network. A candidate turbine $i$ is built if $(i,i+n)$ is selected and not built if $(r,i+n)$ is selected.}
    \label{fig:windfarm-transformation}
\end{figure}

\paragraph{Sets and parameters}
Let $G=(V,A)$ be a directed graph with node set $V=\{1,\dots,2n+1\}$. Node $r=1$ is the substation. The candidate turbine locations are $P=\{2,\ldots,n+1\}$. For each $i\in P$, we introduce a corresponding dummy terminal $i+n$. The set of dummy terminals is $D=\{n+2,\ldots,2n+1\}$, and the set of terminals is $T=\{r\}\cup D$.

The arc set $A$ contains two types of arcs. First, the \emph{cable arcs} form a complete digraph on $\{r\}\cup P$ and represent possible physical cable connections between the substation and candidate turbine locations, or between two candidate locations. The cost of a cable arc combines the cable cost, based on its Euclidean length, with the location-dependent installation cost of the turbine at its head node, reflecting the foundation costs in the underlying data.

Second, for each $i\in P$, we add the auxiliary arcs $(r,i+n)$ and
$(i,i+n)$ to the dummy terminal $i+n$. These arcs have zero cost and are
used only to encode whether turbine $i$ is built: selecting $(i,i+n)$
represents building turbine $i$, whereas selecting $(r,i+n)$ represents
leaving location $i$ empty. 

Just as in SSTP, we let $\delta^-(U)$ denote the arcs entering $U$ for a set $U\subseteq V$. The distance between turbine locations $i$ and $j$ is $d_{ij}$, and $\Delta$ is the minimum allowed distance between two built turbines. The contracted power quota is $Q$, and $\rho$ is the penalty per unit of shortfall.

Uncertainty is represented by a finite scenario set $\mathcal{S}$. Each scenario $s\in\mathcal{S}$ occurs with probability $p^{(s)}$, where $\sum_{s\in\mathcal{S}}p^{(s)}=1$. In scenario $s$, $\xi_i^{(s)}$ is the power yield of turbine $i$. If turbines $i$ and $j$ are both built, $\xi_{ij}^{(s)}$ is the power loss caused by their wake interaction, which is subtracted from their combined production.

\paragraph{Decision variables}
The first-stage binary variable $x_a$ indicates whether arc $a\in A$ is installed. In particular, $x_{(i,i+n)}=1$ means that turbine $i$ is built, while $x_{(r,i+n)}=1$ means that location $i$ is left empty.

For $i,j\in P$ with $i<j$, the continuous variable $m_{ij}\ge 0$ represents the product $x_{(i,i+n)}x_{(j,j')}$. It therefore indicates whether both turbines $i$ and $j$ are built. For each scenario $s\in\mathcal{S}$, the second-stage variable $\sigma^{(s)}\ge 0$ represents the power shortfall relative to the quota.

\paragraph{Formulation}
\begin{subequations}
\begin{align}
\min_{x,\,m,\,\sigma}\quad
  & \sum_{a\in A} c_a\, x_a
  \;+\; \sum_{s \in \mathcal{S}} p^{(s)}\,\rho\,\sigma^{(s)} \label{eq:wf-obj}\\[4pt]
\text{s.t.}\quad
  & \sum_{a\in\delta^-(U)} x_a \ge 1, && \forall i\in P,\ \forall U\subseteq V\setminus\{r\} \text{ with } i+n\in U, \label{eq:wf-conn}\\[2pt]
  & x_{(r,i+n)} + x_{(i,i+n)} = 1, && \forall i\in P, \label{eq:wf-assign}\\[2pt]
  & x_{(i,i+n)} + x_{(j,j+n)} \le 1, && \forall i,j\in P,\ i\neq j,\ d_{ij}<\Delta, \label{eq:wf-prox}\\[2pt]
  & m_{ij} \ge x_{(i,i+n)} + x_{(j,j+n)} - 1, && \forall i,j\in P,\ i<j, \label{eq:wf-mccormick}\\[2pt]
  & \sigma^{(s)} + \sum_{i\in P} \xi_i^{(s)} x_{(i,i+n)} - \sum_{\substack{i,j\in P\\i<j}} \xi_{ij}^{(s)} m_{ij} \ge Q, && \forall s\in\mathcal{S}, \label{eq:wf-demand}\\[2pt]
  & x_a \in \{0,1\}, && \forall\, a\in A, \label{eq:wf-bin}\\[2pt]
  & m_{ij} \ge 0, && \forall\, i,j\in P,\ i<j, \label{eq:wf-mpos}\\[2pt]
  & \sigma^{(s)} \ge 0, && \forall\, s \in \mathcal{S}. \label{eq:wf-spos}
\end{align}
\end{subequations}

The objective~\eqref{eq:wf-obj} minimizes the first-stage cable costs and the expected penalty for power shortfall. The connectivity constraints~\eqref{eq:wf-conn} ensure that every dummy terminal can be reached from the substation. Since these constraints form an exponential family, they are added dynamically during branch-and-cut. Together with~\eqref{eq:wf-assign}, these constraints ensure that a built turbine must be connected to the substation through cable arcs. As for the SSTP, we also add the corresponding symmetric back-cut inequalities to strengthen the LP relaxation.

Constraints~\eqref{eq:wf-mccormick} and~\eqref{eq:wf-mpos} impose the lower McCormick bounds for $x_{(i,i+n)}x_{(j,j+n)}$. The upper McCormick bounds are not required to preserve the optimal objective value. Rearranging~\eqref{eq:wf-demand} gives
\[
    \sigma^{(s)} \geq Q - \sum_{i\in P} \xi_i^{(s)} x_{(i,i+n)} + \sum_{\substack{i,j\in P\\ i<j}} \xi_{ij}^{(s)} m_{ij}.
\]
Since $m_{ij}$ has a nonnegative coefficient in this lower bound and does not appear in the objective, any feasible solution remains feasible after decreasing $m_{ij}$ to its smallest feasible value. Hence, an optimal solution can always be chosen such that
\[
    m_{ij} = \max\{0, x_{(i,i+n)} + x_{(j,j+n)} - 1\} = x_{(i,i+n)} x_{(j,j+n)}.
\]

Finally, constraint~\eqref{eq:wf-demand} requires the produced power, after subtracting wake interference and adding any shortfall $\sigma^{(s)}$, to meet the quota $Q$ in every scenario.

\subsection{Experiments} \label{sec:experiments2}
We use the site~A windfarm layout from the synthetic offshore windfarm dataset of Cazzaro and Pisinger~\citep{Cazzaro2021}, with the \texttt{TNW} wind regime and the $15$--$240$ turbine power curve. The full scenario set contains $|\mathcal{S}|=216$ wind conditions. For $|\mathcal{S}|\in\{50,100,150\}$, we use the first $|\mathcal{S}|$ rows of the wind-condition dataset in their original order. We vary the number of candidate turbine locations over $\{25,50,75,100,125,150\}$.

For this problem, we only consider TULIP-R. The second-stage blocks contain no cuts, while the first-stage block alone does not capture the recourse penalty and therefore provides no useful surrogate. Through the scenario-reduction method from~\cite{feng2016solution}, we reduce the full scenario set to a single representative scenario, $|\mathcal{S}'|=1$. The remaining experimental settings are described in Section~\ref{sec:experimental_design}.

\subsection{Results} \label{sec:case2_results}

\begin{table}[!htbp]
\centering
\caption{Windfarm results by scenario count. For the number of instances solved,
the larger value is shown in bold, with ties highlighted for both methods.
For final optimality gaps, the smaller value is shown in bold.}
\label{tab:windfarm-default}
\footnotesize
\begin{tabular}{rrrrrrcrrr}
\toprule
&& \multicolumn{2}{c}{Solved to optimality}
& \multicolumn{3}{c}{Solved by all methods}
& \multicolumn{3}{c}{Solved by no method} \\
\cmidrule(lr){3-4}
\cmidrule(lr){5-7}
\cmidrule(lr){8-10}
$|\mathcal{S}|$
& $|\mathcal{J}|$
& $|\mathcal{J}_{\mathrm{B\&C}}^\star|$
& $|\mathcal{J}_{\mathrm{R}}^\star|$
& $|\mathcal{J}^{\cap}|$
& $\bar{T}_{\mathrm{B\&C}}$ (s)
& $\bar{\eta}_{\mathrm{R}}$
& $|\mathcal{J}^{0}|$
& $\bar{g}_{\mathrm{B\&C}}$ (\%)
& $\bar{g}_{\mathrm{R}}$ (\%) \\
\midrule
50
& 3
& \textbf{3\,(100\%)}
& \textbf{3\,(100\%)}
& 3 & 261 & 0.44
& 0 & -- & -- \\

100
& 6
& \textbf{3\,(50\%)}
& \textbf{3\,(50\%)}
& 3 & 387 & 0.66
& 3 & \textbf{51.6} & 62.9 \\

150
& 3
& \textbf{3\,(100\%)}
& \textbf{3\,(100\%)}
& 3 & 458 & 0.94
& 0 & -- & -- \\

216 (full)
& 6
& 3\,(50\%)
& \textbf{4\,(67\%)}
& 3 & 238 & 0.72
& 2 & 73.1 & \textbf{46.7} \\
\bottomrule
\end{tabular}
\end{table}

\section{Case study III: stochastic capacitated vehicle routing problem (SCVRP)} \label{sec:scvrp}
\subsection{Problem}
The SCVRP~\citep[Chapter 1.5]{birge_introduction_2011} is a stochastic variant of the vehicle routing problem. For an overview of vehicle routing problems, see~\cite{braekers_vehicle_2016}, and for a review of the SCVRP, see~\cite{berhan_stochastic_2014}. An integer L-shaped algorithm for the SCVRP is presented in~\cite{laporte_integer_2002}. We use the ILP formulation from~\cite{roberti_models_2012}.

A practical example is distribution planning, where a delivery route must be planned before the exact customer demands are known. In the SCVRP considered here, a truck with capacity $C$ must visit the cities $V=\{0,1,\ldots,n\}$, where city $0$ is the depot. We denote the customer cities by $V^*=V\setminus\{0\}$. Traveling from city $i$ to city $j$ has distance $d_{ij}>0$. Customer demands are uncertain and are represented by a finite scenario set $\mathcal{S}$. Scenario $s\in\mathcal{S}$ occurs with probability $p^{(s)}$, and $b_i^{(s)}>0$ denotes the demand of customer $i$ in that scenario. In the first stage, the truck route is chosen before the demands are known. After the demand vector $b^{(s)}$ is revealed, the planned route can be split into additional trips to the depot to satisfy the capacity constraint. Direct travel between two customers is only allowed if that arc was selected in the first-stage route. Demands are served without splitting: every customer is visited exactly once in every scenario. The objective is to minimize the expected total distance traveled while respecting the truck capacity.

\subsection{Example}
To illustrate the SCVRP, we use a small example from \cite[Chapter 1.5]{birge_introduction_2011}. A truck with capacity $C=10$ starts at the depot, visits four customers, and returns to the depot. Figure~\ref{fig:toy_example_graph_scvrp} shows the corresponding complete graph, where the truck can travel directly between any two locations. The demands of customers $1$, $2$, and $4$ are known and equal to $2$ units each. The demand of customer $3$ is uncertain. It is either $1$ or $7$ units, with both outcomes occurring with probability $\tfrac{1}{2}$. The distances between the locations are given in \cite[Table 7]{birge_introduction_2011}.

The optimal first-stage route visits customer $3$ first. If its demand is $1$, the route continues as shown in Figure~\ref{fig:toy_example_scvrp_sol1}. If its demand is $7$, an additional trip to the depot is required, resulting in the route shown in Figure~\ref{fig:toy_example_scvrp_sol2}.

\begin{figure}[!htbp]
    \centering
        \begin{minipage}{0.45\linewidth}
            \centering
             \begin{tikzpicture}[node distance={25mm}, thick, square/.style={draw, regular polygon, regular polygon sides=4}]
                \node[circle, draw, fill=black, minimum size=2mm, inner sep=0pt] (Ep) at (0,0.25) {};
                \node[circle, draw, fill=black, minimum size=2mm, inner sep=0pt] (Ap) at (-1,1) {};
                \node[circle, draw, fill=black, minimum size=2mm, inner sep=0pt] (Bp) at (-0.75,2) {};
                \node[circle, draw, fill=black, minimum size=2mm, inner sep=0pt] (Dp) at (1,1) {};
                \node[circle, draw, fill=black, minimum size=2mm, inner sep=0pt] (Cp) at (0.75,2) {};
                
                \node (E)  [below=0.1cm of Ep] {Depot};
                \node (A) [left=0.1cm of Ap] {\textbf{1} (2)};
                \node (B) [left=0.1cm of Bp] {\textbf{2} (2)};
                \node (D) [right=0.1cm of Dp] {\textbf{4} (2)};
                \node (C) [right=0.1cm of Cp] {\textbf{3} (1 or 7)};
                \node (extra) [above=0.03cm of C] {$p^{(1)}=p^{(2)}=\frac{1}{2}$};

                \draw (Ap) -- (Bp);
                \draw (Ap) -- (Cp);
                \draw (Ap) -- (Dp);
                \draw (Ap) -- (Ep);
                \draw (Bp) -- (Cp);
                \draw (Bp) -- (Dp);
                \draw (Bp) -- (Ep);
                \draw (Cp) -- (Dp);
                \draw (Cp) -- (Ep);
                \draw (Dp) -- (Ep);
            \end{tikzpicture}
            \subcaption{Problem instance.}
            \label{fig:toy_example_graph_scvrp}
        \end{minipage}
        \vfill
        \begin{minipage}{0.45\linewidth}
            \centering
             \begin{tikzpicture}[node distance={25mm}, thick, square/.style={draw, regular polygon, regular polygon sides=4}]
                \node[circle, draw, fill=black, minimum size=2mm, inner sep=0pt] (Ep) at (0,0.25) {};
                \node[circle, draw, fill=black, minimum size=2mm, inner sep=0pt] (Ap) at (-1,1) {};
                \node[circle, draw, fill=black, minimum size=2mm, inner sep=0pt] (Bp) at (-0.75,2) {};
                \node[circle, draw, fill=black, minimum size=2mm, inner sep=0pt] (Dp) at (1,1) {};
                \node[circle, draw, fill=black, minimum size=2mm, inner sep=0pt] (Cp) at (0.75,2) {};

                \draw[->] (Ep) -- (Cp);
                \draw (Ap) -- (Bp);
                \draw (Ap) -- (Dp);
                \draw (Bp) -- (Cp);
                \draw (Dp) -- (Ep);
            \end{tikzpicture}
            \subcaption{Optimal route for scenario $1$.}
            \label{fig:toy_example_scvrp_sol1}
        \end{minipage}
        \hfill
        \begin{minipage}{0.45\linewidth}
            \centering
             \begin{tikzpicture}[node distance={25mm}, thick, square/.style={draw, regular polygon, regular polygon sides=4}]
                \node[circle, draw, fill=black, minimum size=2mm, inner sep=0pt] (Ep) at (0,0.25) {};
                \node[circle, draw, fill=black, minimum size=2mm, inner sep=0pt] (Ap) at (-1,1) {};
                \node[circle, draw, fill=black, minimum size=2mm, inner sep=0pt] (Bp) at (-0.75,2) {};
                \node[circle, draw, fill=black, minimum size=2mm, inner sep=0pt] (Dp) at (1,1) {};
                \node[circle, draw, fill=black, minimum size=2mm, inner sep=0pt] (Cp) at (0.75,2) {};

                \draw[->] (Ep) -- (Cp);
                \draw[->] (Ep) -- (Ap);
                \draw (Ap) -- (Dp);
                \draw (Dp) -- (Ep);
                \draw (Bp) -- (Cp);
                \draw (Bp) -- (Ep);
            \end{tikzpicture}
            \subcaption{Optimal route for scenario $2$.}
            \label{fig:toy_example_scvrp_sol2}
        \end{minipage}
        \caption{Small SCVRP example from \cite[Chapter 1.5]{birge_introduction_2011}. Bold numbers denote the customers, and numbers in parentheses denote their demands. Customer $3$ has demand $1$ or $7$, with both scenarios occurring with probability $\tfrac{1}{2}$.}
    \label{fig:toy_example_scvrp}
\end{figure}

\subsection{Model}
\paragraph{Sets and parameters}
The truck operates on the node set $V=\{0,1,\dots,n\}$, where node $0$ is the depot and $V^{*}=V\setminus\{0\}$ is the set of customers. Let $A=\{(i,j)\in V\times V:i\neq j\}$ be the set of directed arcs, excluding self-loops. Traveling directly from node $i$ to node $j$ has distance $d_{ij}>0$, and the truck has capacity $C$. Uncertainty is represented by a finite scenario set $\mathcal{S}$. Each scenario $s\in\mathcal{S}$ occurs with probability $p^{(s)}$. In scenario $s$, customer $i$ has demand $b_i^{(s)}>0$.

\paragraph{Decision variables}
The first-stage variables $x_{ij}\in\{0,1\}$ describe the planned route. Variable $x_{ij}$ equals $1$ if the truck is planned to travel directly from node $i$ to node $j$. Because the depot degrees are bounded below rather than fixed at one, $x$ may describe one or several depot-rooted routes.

After the demands are revealed, the second-stage variables $y_{ij}^{(s)}\in\{0,1\}$ describe the route followed in scenario $s$. Variable $y_{ij}^{(s)}$ equals $1$ if the truck travels directly from node $i$ to node $j$ in that scenario. The second-stage route may include additional trips to the depot.

\paragraph{Formulation}
\begin{subequations} \label{eq:scvrp}
\begin{alignat}{3}
  \min \quad
    & \sum_{s \in \mathcal{S}} p^{(s)} \sum_{(i,j)\in A} d_{ij}\, y_{ij}^{(s)}
    && \label{eq:scvrp-obj}\\
  \text{s.t.}\quad
    & \sum_{j \in V,\, j \neq i} x_{ij} = 1, && \forall\, i \in V^{*}, \label{eq:scvrp-out}\\
    & \sum_{i \in V,\, i \neq j} x_{ij} = 1, && \forall\, j \in V^{*}, \label{eq:scvrp-in}\\
    & \sum_{j \in V^{*}} x_{0j} \geq 1, \qquad \sum_{i \in V^{*}} x_{i0} \geq 1, && \label{eq:scvrp-depot}\\
    & \sum_{j \in V,\, j \neq i} y_{ij}^{(s)} = 1, && \forall\, i \in V^{*},\ \forall\, s \in \mathcal{S}, \label{eq:scvrp-y-out}\\
    & \sum_{i \in V,\, i \neq j} y_{ij}^{(s)} = 1, && \forall\, j \in V^{*},\ \forall\, s \in \mathcal{S}, \label{eq:scvrp-y-in}\\
    & \sum_{j \in V^{*}} y_{0j}^{(s)} \geq 1, \qquad \sum_{i \in V^{*}} y_{i0}^{(s)} \geq 1, && \forall\, s \in \mathcal{S}, \label{eq:scvrp-y-depot}\\
    & \sum_{i \in Q} \sum_{j \in Q} x_{ij} \leq |Q| - 1, && \forall\, Q \subseteq V^{*},\ |Q| \geq 2, \label{eq:scvrp-subtour}\\
    & \sum_{i \in Q} \sum_{j \in V \setminus Q} \bigl(y^{(s)}_{ij} + y^{(s)}_{ji}\bigr) \geq 2 \left\lceil \frac{\sum_{i \in Q} b_i^{(s)}}{C} \right\rceil, && \forall\, Q \subseteq V^{*},\ Q \neq \emptyset,\ \forall\, s \in \mathcal{S}, \label{eq:scvrp-cap}\\
    & y_{ij}^{(s)} \leq x_{ij}, && \forall\, i,j \in V^{*},\ \forall\, s \in \mathcal{S}, \label{eq:scvrp-couple}\\
    & x_{ij},\, y_{ij}^{(s)} \in \{0,1\}, && \forall\, (i,j)\in A,\ \forall\, s \in \mathcal{S}. \label{eq:scvrp-bin}
\end{alignat}
\end{subequations}

The objective~\eqref{eq:scvrp-obj} minimizes the expected distance traveled. Constraints~\eqref{eq:scvrp-out} and~\eqref{eq:scvrp-in} require every customer to be entered and left exactly once in the planned route. Constraints~\eqref{eq:scvrp-depot} require the route to use the depot. Likewise, constraints~\eqref{eq:scvrp-y-out} and~\eqref{eq:scvrp-y-in} impose one outgoing and one incoming realized arc at every customer in every scenario, and constraints~\eqref{eq:scvrp-y-depot} require every realized solution to use the depot. Together with the exclusion of self-loops, these constraints enforce single-visit, unsplit service. The subtour-elimination constraints~\eqref{eq:scvrp-subtour} prevent planned routes that are disconnected from the depot~\citep{roberti_models_2012}. Since the latter form an exponential family, they are generated dynamically during branch-and-cut when violated. The capacity constraints~\eqref{eq:scvrp-cap} are rounded-capacity cuts~\citep[Chapter 3]{toth_vehicle_2002,hernandez-perez_branch-and-cut_2004}. For every nonempty subset of customers, there must be enough arcs to cross the subset's boundary to serve its total demand without exceeding the truck capacity. These constraints also form an exponential family and are added dynamically during branch-and-cut. Finally, the coupling constraints~\eqref{eq:scvrp-couple} link the second-stage route to the planned first-stage route. Between customers, the second-stage route may only use planned arcs. Recourse, therefore, consists of adding additional trips to the depot as needed.

\subsection{Experiments} \label{sec:experiments3}
We use three capacitated vehicle-routing instances: eil7, eil13, and eil23. The instances eil13 and eil23 correspond to the TSPLIB/CVRP instances E-n13-k4 and E-n23-k3~\cite{reinelt_gerhard_tsplib_2024}. The smaller eil7 instance is constructed in the same style, with unit nominal demands. Table~\ref{tab:ctsp_instance} summarizes the instances.

Travel distances are fixed across scenarios, while customer demands are uncertain. Let $B_i$ denote the demand of customer $i$ in the original CVRP instance, and let $C_0$ denote the vehicle capacity in that instance. Following~\cite{juan_using_2011}, for every customer $i\in V^*$ and scenario $s\in\mathcal{S}$ we draw independently
\begin{equation}
Z_i^{(s)} \sim \mathcal{N}\!\left(\frac{B_i}{C_0},\left(0.25\frac{B_i}{C_0}\right)^2\right),
\qquad
\widetilde b_i^{(s)}=\exp\!\left(Z_i^{(s)}\right).
\label{eq:scvrp-demand-draw}
\end{equation}
The experimental vehicle capacity is then calibrated from the untruncated finite scenario sample to obtain three vehicle loads on average, after which the single-customer cap is applied:
\begin{equation}
C=\frac{1}{3|\mathcal{S}|}\sum_{s\in\mathcal{S}}\sum_{i\in V^*}\widetilde b_i^{(s)},
\qquad
b_i^{(s)}=\min\!\left\{\widetilde b_i^{(s)},0.95C\right\}.
\label{eq:scvrp-demand-calibration}
\end{equation}
No rounding is applied. The cap guarantees that each customer can be served in a single visit; it was inactive in all demand sets used in the experiments (the largest realized customer demand was less than $0.66C$). Consequently, the sample-average total demand equals exactly three vehicle loads in every tested instance. This calibration preserves a meaningful multi-route stochastic routing problem; choosing a capacity that permits all expected demand to be served in one route would remove much of the stochastic routing structure that we aim to study.

\begin{table}[!htbp]
    \centering
    \caption{SCVRP instances used in the experiments. Here, $B_i$ denotes the nominal demand of customer $i$.}
    \label{tab:ctsp_instance}
    \begin{tabular}{@{}lccc@{}}
    \toprule
    Instance & $|V|$ & $\min(B_i)$ & $\max(B_i)$ \\
    \midrule
    eil7  &  7 &    1 &    1 \\
    eil13 & 13 & 1100 & 1900 \\
    eil23 & 23 &   60 & 4100 \\
    \bottomrule
    \end{tabular}
\end{table}

For each instance, we consider scenario counts $|\mathcal{S}|\in\{10,20,50\}$ and three independent demand sets. We compare regular branch-and-cut with TULIP-R, TULIP-F, TULIP-S, and TULIP-C. The second-stage surrogates omit the coupling constraints $y_{ij}^{(s)}\le x_{ij}$ and optimize each scenario's routes independently, with their construction and solution times included in $t_{\mathrm{sur}}$. The remaining experimental settings are described in Section~\ref{sec:experimental_design}.

\subsection{Results} \label{sec:case3_results}

\begin{table}[H]
\centering
\caption{Final optimality gaps for SCVRP instances not solved by any method within the time limit. Bold values indicate the best performance within each scenario count, and ties are all highlighted.}
\label{tab:gapdiff-truck-default}
\small
\begin{tabular}{crccccc}
\toprule
$|\mathcal{S}|$ & $|\mathcal{J}^0|$ & $\bar{g}_{\mathrm{B\&C}}$ (\%) & $\overline{\Delta g}_\mathrm{R}$ & $\overline{\Delta g}_\mathrm{F}$ & $\overline{\Delta g}_\mathrm{S}$ & $\overline{\Delta g}_\mathrm{C}$ \\
\midrule
$10$ & 4 & 2.1 & \textbf{+0.4} & \textbf{+0.4} & +0.0 & \textbf{+0.4} \\
$20$ & 6 & 5.1 & \textbf{+0.8} & +0.6 & 0.0 & +0.6 \\
$50$ & 6 & 23.7 & +9.9 & \textbf{+17.4} & 0.0 & +17.3 \\
\bottomrule
\end{tabular}
\end{table}

\begin{table}[H]
\centering
\caption{Number of SCVRP instances solved and mean speedup by scenario count. Bold values indicate the best performance within each scenario count, and ties are all highlighted.}
\label{tab:solved-speedup-truck-default}
\footnotesize
\setlength{\tabcolsep}{3pt}
\begin{tabular}{rrrcccccccccc}
\toprule
&&& \multicolumn{2}{c}{B\&C} & \multicolumn{2}{c}{TULIP-R} & \multicolumn{2}{c}{TULIP-F} & \multicolumn{2}{c}{TULIP-S} & \multicolumn{2}{c}{TULIP-C} \\
\cmidrule(lr){4-5}\cmidrule(lr){6-7}\cmidrule(lr){8-9}\cmidrule(lr){10-11}\cmidrule(lr){12-13}
$|\mathcal{S}|$ & $|\mathcal{J}|$ & $|\mathcal{J}^{\cap}|$ & $|\mathcal{J}_{\mathrm{B\&C}}^\star|$ & $\bar{T}_{\mathrm{B\&C}}$ (s) & $|\mathcal{J}_\mathrm{R}^\star|$ & $\bar{\eta}_{\mathrm{R}}$ & $|\mathcal{J}_\mathrm{F}^\star|$ & $\bar{\eta}_\mathrm{F}$ & $|\mathcal{J}_\mathrm{S}^\star|$ & $\bar{\eta}_\mathrm{S}$ & $|\mathcal{J}_\mathrm{C}^\star|$ & $\bar{\eta}_\mathrm{C}$ \\
\midrule
$10$ & 9 & 4 &
4\,(44\%) & 159 &
4\,(44\%) & 0.28 &
4\,(44\%) & 0.31 &
\textbf{5\,(56\%)} & 0.89 &
4\,(44\%) & 0.28 \\

$20$ & 9 & 3 &
\textbf{3\,(33\%)} & 0.7 &
\textbf{3\,(33\%)} & 0.06 &
\textbf{3\,(33\%)} & 0.39 &
\textbf{3\,(33\%)} & 0.63 &
\textbf{3\,(33\%)} & 0.33 \\

$50$ & 9 & 3 &
\textbf{3\,(33\%)} & 2.3 &
\textbf{3\,(33\%)} & 0.15 &
\textbf{3\,(33\%)} & 0.64 &
\textbf{3\,(33\%)} & 0.77 &
\textbf{3\,(33\%)} & 0.56 \\
\bottomrule
\end{tabular}
\end{table}

\section{Case study IV: storage, layout, and routing problem (SLR)} \label{sec:slr}
\subsection{Problem}
We consider the storage, layout, and routing problem introduced in~\citep{narum_problem-based_2024}. A set of products must be assigned to storage locations in a warehouse before customer orders are known. Uncertainty is represented by a finite set of scenarios, where each scenario specifies a batch of products that must be collected. In the first stage, products are assigned to storage locations. After a scenario is revealed, a picking route is constructed to visit the locations containing the requested products. The objective is to choose the product locations such that the expected picking time over all scenarios is minimized.

We consider two versions of the problem. In the base version, denoted by SLR, the first stage only assigns products to storage locations. In \slrp, the first stage additionally installs a connected backbone between the stocked locations, which can be used by the second-stage picking routes without an additional penalty. Routes may also use connections outside the backbone, but at an additional cost.

\subsection{Example SLR}
To illustrate the problem, consider a warehouse with three storage locations arranged in a $3\times1$ grid. The first location contains the door, where every picking route starts. Following~\cite{narum_problem-based_2024}, item $1$ represents the door and is included in every batch, ensuring that every picking route starts and ends at the door. This is a modeling choice to ensure that every scenario requires a visit to the door, from which the picking route starts. Travel times between locations are given by the Manhattan distance. There are five items in total, and each location can store at most two items. We consider two equally likely scenarios. The first batch contains items $1$, $3$, and $5$, while the second contains items $1$, $2$, $3$, and $4$.

One optimal layout stores items $1$ and $3$ in the first location,
items $2$ and $4$ in the second, and item $5$ in the third. In the first scenario, the requested items are stored in the first and third locations, giving a picking-route cost of $2$. In the second scenario, the requested items are stored in the first and second locations, giving a picking-route cost of $1$. The corresponding expected objective value is therefore $1.5$.

\subsection{Model}
\paragraph{Sets and parameters}
Let $L$ be the set of storage locations and $I$ the set of items. Location $1$ is the facility door, and item $1$ is fixed at this location. We write $L^*=L\setminus\{1\}$. Let
\[
    J=\{(l_1,l_2):l_1,l_2\in L,\ l_1<l_2\}
\]
be the set of unordered pairs of locations. The travel time between the two locations in pair $j\in J$ is $d_j$. Let $\bar c$ denote the maximum number of items that can be stored at one location. Uncertainty is represented by a finite scenario set $\mathcal{S}$, where scenario $s$ occurs with probability $p^{(s)}$. The binary parameter $\xi_i^{(s)}$ equals $1$ if item $i$ is part of the batch in scenario $s$, and $0$ otherwise. Finally, $M$ is a sufficiently large constant.

\paragraph{Decision variables}
The first-stage binary variable $w_{il}$ equals $1$ if item $i$ is stored at location $l$. Duplicate sampled scenarios are retained as separate recourse blocks, so these experiments study computational scaling with sample size rather than with the number of distinct demand outcomes. Then, for each scenario $s\in\mathcal{S}$, the second-stage binary variable $z_j^{(s)}$ equals $1$ if the picking route uses the edge between the two locations in pair $j$. Variable $v_l^{(s)}$ indicates whether location $l$ contains at least one item required in scenario $s$. Finally, $\tau^{(s)}$ equals $1$ when the requested items are stored across more than two locations, in which case the picking route is modeled as a closed tour.

\paragraph{Formulation}
\begingroup
\setlength{\jot}{1pt}
\begin{subequations} \label{eq:slr}
    \begin{alignat}{3}
        \min_{w,z,v,\tau} \quad & \sum_{s \in \mathcal{S}} p^{(s)} \sum_{j \in J} d_j z^{(s)}_j & \label{eq:slr_1}\\
        \mbox{s.t.} \quad & \sum_{i \in I} w_{il} \leq \bar{c} \qquad & \forall l \in L \label{eq:slr_2}\\
        & \sum_{l \in L} w_{il} = 1 & \forall i \in I \label{eq:slr_3}\\
        & w_{11} = 1 & \label{eq:slr_4}\\
        & \sum_{l \in L} v^{(s)}_l - M\,\tau^{(s)} \leq 2 \qquad & \forall s \in \mathcal{S} \label{eq:slr_5}\\
        & v^{(s)}_l \geq \frac{\sum_{i \in I} w_{il}\,\xi^{(s)}_i}{\bar{c}} & \forall l \in L,\ \forall s \in \mathcal{S} \label{eq:slr_6}\\
        & \sum_{\substack{j \in J:\ l \in j}} z^{(s)}_{j} \geq v^{(s)}_l & \forall l \in L^*,\ \forall s \in \mathcal{S} \label{eq:slr_7}\\
        & \sum_{\substack{j \in J:\ 1 \in j}} z^{(s)}_{j} \geq \frac{\sum_{l \in L} v^{(s)}_l - 1}{|I|} & \forall s \in \mathcal{S} \label{eq:slr_8}\\
        & \sum_{j \in J} z^{(s)}_j = \sum_{l \in L} v^{(s)}_l - 1 + \tau^{(s)} \quad & \forall s \in \mathcal{S} \label{eq:slr_9}\\
        & \sum_{\substack{j \in J:\ l \in j}} z^{(s)}_{j} \leq 2\,v^{(s)}_l \quad & \forall l \in L,\ \forall s \in \mathcal{S} \label{eq:slr_10}\\
        & \sum_{\substack{j \in J:\ j \subseteq L'}} z^{(s)}_j \leq | L' | - 1 \quad & \forall L' \subseteq L^* : |L'| \geq 2,\ \forall s \in \mathcal{S} \label{eq:slr_15}\\
        & w_{il} \in \{0,1\} & \forall i \in I,\ \forall l \in L \label{eq:slr_11}\\
        & z^{(s)}_j \in \{0,1\} & \forall j \in J,\ \forall s \in \mathcal{S} \label{eq:slr_12}\\
        & v^{(s)}_l \in \{0,1\} & \forall l \in L,\ \forall s \in \mathcal{S} \label{eq:slr_13}\\
        & \tau^{(s)} \in \{0,1\} & \forall s \in \mathcal{S} \label{eq:slr_14}
    \end{alignat}
\end{subequations}
\endgroup

The objective~\eqref{eq:slr_1} minimizes the expected travel time of the picking route. Constraint~\eqref{eq:slr_2} limits the number of items stored at each location, while~\eqref{eq:slr_3} ensures that every item is assigned to exactly one location. Following~\cite{narum_problem-based_2024}, item $1$ is included in every batch. Constraint~\eqref{eq:slr_4} fixes this item at the door, thereby ensuring that every picking route includes the door. Constraint~\eqref{eq:slr_6} determines which locations are active in each scenario. If location $l$ contains at least one item required in scenario $s$, then $v_l^{(s)}$ must equal $1$. Constraints~\eqref{eq:slr_7} and~\eqref{eq:slr_8} ensure that the route uses the active locations and leaves the door whenever more than one location is active. In Constraint~\eqref{eq:slr_8}, the right-hand side is zero when only the door is required and positive when at least one additional location must be visited. Since the left-hand side is an integer, any positive right-hand side forces at least one route edge incident to the door. Constraint~\eqref{eq:slr_5} sets $\tau^{(s)}=1$ when the requested items are stored across more than two locations. Together with Constraint~\eqref{eq:slr_9}, this determines the number of edges in the picking route. Following~\cite{narum_problem-based_2024}, when only the door and one other location are active, the route traverses the same undirected connection in both directions. Since the binary variable $z_j^{(s)}$ cannot represent using the same edge twice, Constraint~\eqref{eq:slr_9} selects a single edge between them, and this case is treated separately. If three or more locations are active, the selected edges instead form a closed tour that starts and ends at the door. Constraint~\eqref{eq:slr_10} limits each active location to at most two incident route edges. Finally, Constraint~\eqref{eq:slr_15} eliminates subtours. They form an exponential family and are therefore generated dynamically during branch-and-cut. Constraints~\eqref{eq:slr_11}--\eqref{eq:slr_14} state the integrality requirements.

\subsection{First-stage routing variant}
In the routing version, the first stage additionally installs a backbone network,
representing infrastructure such as a conveyor system or a set of main aisles.
This infrastructure is installed before future orders are known and can be used by all second-stage picking routes. A picking route may also use a connection outside the backbone, but incurs an additional penalty for doing so. We denote this version by \slrp{} and the base assignment model~\eqref{eq:slr} by SLR. The routing version introduces an exponential family of first-stage connectivity constraints, making TULIP-F and TULIP-C applicable.

We introduce two first-stage binary variables. Variable $u_k$ equals $1$ if location $k$ contains at least one item, and $w_a^1$ equals $1$ if backbone arc $a$ is installed. In the second stage, $y_j^{(s)}$ equals $1$ if picking edge $j$ is used outside the installed backbone in scenario $s$. The parameter $\kappa$ is the backbone installation cost per unit of distance, and $\rho$ is the penalty for using an edge outside the backbone. To model the connected backbone, we use a directed representation rooted at the door. Let $A=\{(i,j):i,j\in L,\ i\neq j\}$ be the set of directed arcs between locations, and let $\delta^-(S)$ denote the set of arcs entering $S\subseteq L$. The arc directions are only a modeling convention: a backbone connection can be used by
a second-stage picking route in either direction. Let $r=1$ denote the root at the door. The base model~\eqref{eq:slr} is then augmented with
\begin{subequations} \label{eq:slr_route}
    \begin{alignat}{3}
        & u_k \leq \sum_{i \in I} w_{ik} \qquad & \forall k \in L \label{eq:slrr_1}\\
        & \bar{c}\,u_k \geq \sum_{i \in I} w_{ik} & \forall k \in L \label{eq:slrr_2}\\
        & \sum_{a \in \delta^-(k)} w^{1}_a = u_k & \forall k \in L \setminus \{r\} \label{eq:slrr_3}\\
        & \sum_{a \in \delta^-(r)} w^{1}_a = 0 & \label{eq:slrr_4}\\
        & \sum_{a \in \delta^-(S)} w^{1}_a \geq u_t \quad & \forall S \subseteq L \setminus \{r\},\ \forall t \in S \label{eq:slrr_5}\\
        & z^{(s)}_{(i,j)} \leq w^{1}_{(i,j)} + w^{1}_{(j,i)} + y^{(s)}_{(i,j)} \quad & \forall (i,j) \in J,\ \forall s \in \mathcal{S} \label{eq:slrr_6}\\
        & u_k,\ w^{1}_a \in \{0,1\} & \forall k \in L,\ \forall a \in A \label{eq:slrr_7}\\
        & y^{(s)}_j \in \{0,1\} & \forall j \in J,\ \forall s \in \mathcal{S}. \label{eq:slrr_8}
    \end{alignat}
\end{subequations}
and the objective~\eqref{eq:slr_1} is replaced by
\begin{equation} \label{eq:slr_route_obj}
    \min \quad \kappa \sum_{a \in A} d_a\, w^{1}_a + \sum_{s \in \mathcal{S}} p^{(s)} \left(\sum_{j \in J} d_j z^{(s)}_j + \rho \sum_{j \in J} y^{(s)}_j \right).
\end{equation}

Constraints~\eqref{eq:slrr_1} and~\eqref{eq:slrr_2} determine whether a location is stocked. Constraints~\eqref{eq:slrr_3} and~\eqref{eq:slrr_4} assign one incoming backbone arc to every stocked non-root location and none to the root. The connectivity cuts~\eqref{eq:slrr_5} ensure that every stocked location is connected to the root. These constraints form an exponential family and are generated dynamically during branch-and-cut. Constraint~\eqref{eq:slrr_6} couples the two stages. A picking edge can be used if the corresponding backbone connection is installed in either direction. Otherwise, the edge can still be used by setting $y_j^{(s)}=1$ and paying the additional penalty $\rho$. The objective~\eqref{eq:slr_route_obj} therefore trades off the cost of installing the backbone against the expected cost of using edges outside it.

\paragraph{Example \slrp}
Consider a $2\times2$ grid with four locations as shown in Figure~\ref{fig:stp+_example}. Location $1=(0,0)$ is the door, and the remaining locations are $2=(1,0)$, $3=(0,1)$, and $4=(1,1)$. Distances are given by the Manhattan distance. There are four items; each location can contain one item, and item $1$ is fixed at the door. We set $\kappa=1$ and $\rho=2$. There are two equally likely scenarios. Scenario~A requires items $\{1,2,3\}$, while scenario~B requires items $\{1,2,4\}$. 

In the optimal solution, item $i$ is stored at location $i$. The first stage installs the backbone
\[
    \{1\!\to\!2,\;1\!\to\!3,\;2\!\to\!4\},
\]
which has total installation cost $3$. In scenario~A, the picking route is $1\!-\!2\!-\!3\!-\!1$. Its travel distance is $4$. The edges $(1,2)$ and $(1,3)$ are available on the backbone, while $(2,3)$ requires a detour and therefore adds the penalty $\rho=2$. In scenario~B, the route is $1\!-\!2\!-\!4\!-\!1$, again with travel distance $4$. Here, $(1,2)$ and $(2,4)$ are available on the backbone, while $(1,4)$ requires a detour. The total expected cost is therefore
\[
    3
    + \tfrac12(4+2)
    + \tfrac12(4+2)
    = 9.
\]

The example illustrates the trade-off introduced by the routing variant. Installing a backbone connection has a first-stage cost, but allows that connection to be used in any scenario without an additional penalty. Not installing it saves this upfront cost, but incurs the penalty $\rho$ whenever the connection is needed by a second-stage picking route. The first-stage decision therefore balances the cost of installing backbone connections against their expected future use.

\begin{figure}[!htbp]
\centering
\begin{tikzpicture}[>=Stealth, scale=1.1,
   loc/.style={circle,draw,fill=white,minimum size=6mm,inner sep=0pt},
   bb/.style={gray!55,line width=2.2pt,->},
   route/.style={blue!70!black,line width=1.3pt},
   detour/.style={blue!70!black,line width=1.3pt,dashed}]
 \node[loc] (1) at (0,0) {1}; \node[loc] (2) at (2,0) {2};
 \node[loc] (3) at (0,2) {3}; \node[loc] (4) at (2,2) {4};
 \draw[bb] (1) -- (2); \draw[bb] (1) -- (3); \draw[bb] (2) -- (4);
 \draw[route] (1) -- (2); \draw[route] (1) -- (3); \draw[detour] (2) -- (3);
 \node at (1,-0.8) {(a) scenario A: tour $1\!-\!2\!-\!3$};
\end{tikzpicture}
\hfill
\begin{tikzpicture}[>=Stealth, scale=1.1,
   loc/.style={circle,draw,fill=white,minimum size=6mm,inner sep=0pt},
   bb/.style={gray!55,line width=2.2pt,->},
   route/.style={blue!70!black,line width=1.3pt},
   detour/.style={blue!70!black,line width=1.3pt,dashed}]
 \node[loc] (1) at (0,0) {1}; \node[loc] (2) at (2,0) {2};
 \node[loc] (3) at (0,2) {3}; \node[loc] (4) at (2,2) {4};
 \draw[bb] (1) -- (2); \draw[bb] (1) -- (3); \draw[bb] (2) -- (4);
 \draw[route] (1) -- (2); \draw[route] (2) -- (4); \draw[detour] (1) -- (4);
 \node at (1,-0.8) {(b) scenario B: tour $1\!-\!2\!-\!4$};
\end{tikzpicture}
\caption{Example of the \slrp\ routing variant on a $2\times2$ grid. Thick gray arcs show the first-stage backbone. Solid blue edges are picking edges that use the backbone, while dashed blue edges are off-backbone edges that incur the detour penalty $\rho$.}
\label{fig:stp+_example}
\end{figure}

\subsection{Experiments} \label{sec:experiments4}
We consider instances with $|I|=8$ items and $|L|=6$ storage locations arranged in a $3\times2$ grid. Each location can store at most $\bar{c}=2$ items. Travel times between locations are given by the Manhattan distance.

Following~\cite{narum_problem-based_2024}, uncertain batches are generated as correlated binary item demands. Item $1$, which is located at the door, is included in every batch. For each of the remaining items, correlated normal draws are generated and thresholded such that the item is included with probability $5/|I|$. The covariance matrix is generated randomly and projected to be positive semidefinite.

For SLR, we consider scenario counts
\[
    |\mathcal{S}|\in\{50,100,150,200,300,500,750,1000\},
\]
while for \slrp{} we consider
\[
    |\mathcal{S}|\in\{50,100,150,200,300,500\}.
\]
For each scenario count, we generate $10$ independent demand instances. For \slrp, we set the backbone installation cost to $\kappa=1$ and the detour penalty to
\[
    \rho=2(\sqrt{|L|}-1).
\]

For SLR, we compare regular branch-and-cut with TULIP-R and TULIP-S. For \slrp, we compare regular branch-and-cut with TULIP-R, TULIP-F, TULIP-S, and TULIP-C. The remaining experimental settings are described in Section~\ref{sec:experimental_design}.

\paragraph{First-stage choice for TULIP-S}
TULIP-S solves individual second-stage blocks and therefore requires the first-stage decisions to be fixed beforehand. For SLR and \slrp, we use the corresponding expected-value solution for these decisions. Each scenario surrogate is then solved with the first-stage decisions fixed, and the second-stage cuts collected during these solves are transferred to the full problem as described in Section~\ref{sec:algorithm}.

\clearpage
\subsection{Results} \label{sec:case4_results}
\begin{table}[!htbp]
\centering
\caption{Number of SLR instances solved and mean speedup by scenario count. Bold values indicate the best performance within each scenario count, and ties are all highlighted.}
\label{tab:solved-speedup-slr-plain}
\begin{tabular}{rrrcccccc}
\toprule
&&& \multicolumn{2}{c}{B\&C} & \multicolumn{2}{c}{TULIP-R} & \multicolumn{2}{c}{TULIP-S} \\
\cmidrule(lr){4-5}\cmidrule(lr){6-7}\cmidrule(lr){8-9}
$|\mathcal{S}|$ & $|\mathcal{J}|$ & $|\mathcal{J}^{\cap}|$ & $|\mathcal{J}_{\mathrm{B\&C}}^\star|$ & $\bar{T}_{\mathrm{B\&C}}$ (s) & $|\mathcal{J}_\mathrm{R}^\star|$ & $\bar{\eta}_{\mathrm{R}}$ & $|\mathcal{J}_\mathrm{S}^\star|$ & $\bar{\eta}_{\mathrm{S}}$ \\
\midrule
$50$ & 10 & 10 & \textbf{10\,(100\%)} & 3.9 & \textbf{10\,(100\%)} & 0.25 & \textbf{10\,(100\%)} & 0.76 \\
$100$ & 10 & 10 & \textbf{10\,(100\%)} & 10 & \textbf{10\,(100\%)} & 0.46 & \textbf{10\,(100\%)} & 0.82 \\
$150$ & 10 & 10 & \textbf{10\,(100\%)} & 22 & \textbf{10\,(100\%)} & 0.64 & \textbf{10\,(100\%)} & 0.87 \\
$200$ & 10 & 10 & \textbf{10\,(100\%)} & 41 & \textbf{10\,(100\%)} & 0.72 & \textbf{10\,(100\%)} & 0.89 \\
$300$ & 10 & 10 & \textbf{10\,(100\%)} & 90 & \textbf{10\,(100\%)} & 0.85 & \textbf{10\,(100\%)} & 0.94 \\
$500$ & 10 & 10 & \textbf{10\,(100\%)} & 408 & \textbf{10\,(100\%)} & 0.95 & \textbf{10\,(100\%)} & 0.97 \\
$750$ & 10 & 2 & \textbf{6\,(60\%)} & 1729 & 2\,(20\%) & \textbf{1.15} & \textbf{6\,(60\%)} & 0.99 \\
$1000$ & 10 & 0 & \textbf{1\,(10\%)} & -- & 0\,(0\%) & -- & \textbf{1\,(10\%)} & -- \\
\bottomrule
\end{tabular}
\end{table}

\begin{table}[!htbp]
\centering
\caption{Number of \slrp\ instances solved and mean speedup by scenario count. Bold values indicate the best performance within each scenario count, and ties are all highlighted.}
\label{tab:solved-speedup-slr-routing}
\resizebox{\linewidth}{!}{
\begin{tabular}{rrrcccccccccc}
\toprule
&&& \multicolumn{2}{c}{B\&C} & \multicolumn{2}{c}{TULIP-R} & \multicolumn{2}{c}{TULIP-F} & \multicolumn{2}{c}{TULIP-S} & \multicolumn{2}{c}{TULIP-C} \\
\cmidrule(lr){4-5}\cmidrule(lr){6-7}\cmidrule(lr){8-9}\cmidrule(lr){10-11}\cmidrule(lr){12-13}
$|\mathcal{S}|$ & $|\mathcal{J}|$ & $|\mathcal{J}^{\cap}|$ & $|\mathcal{J}_{\mathrm{B\&C}}^\star|$ & $\bar{T}_{\mathrm{B\&C}}$ (s) & $|\mathcal{J}_\mathrm{R}^\star|$ & $\bar{\eta}_{\mathrm{R}}$ & $|\mathcal{J}_\mathrm{F}^\star|$ & $\bar{\eta}_{\mathrm{F}}$ & $|\mathcal{J}_\mathrm{S}^\star|$ & $\bar{\eta}_{\mathrm{S}}$ & $|\mathcal{J}_\mathrm{C}^\star|$ & $\bar{\eta}_{\mathrm{C}}$ \\
\midrule
$50$ & 10 & 10 & \textbf{10\,(100\%)} & 41 & \textbf{10\,(100\%)} & 0.72 & \textbf{10\,(100\%)} & 0.93 & \textbf{10\,(100\%)} & 0.94 & \textbf{10\,(100\%)} & 0.88 \\
$100$ & 10 & 10 & \textbf{10\,(100\%)} & 224 & \textbf{10\,(100\%)} & \textbf{1.08} & \textbf{10\,(100\%)} & 0.98 & \textbf{10\,(100\%)} & 0.97 & \textbf{10\,(100\%)} & 0.95 \\
$150$ & 10 & 10 & \textbf{10\,(100\%)} & 380 & \textbf{10\,(100\%)} & 0.95 & \textbf{10\,(100\%)} & 0.99 & \textbf{10\,(100\%)} & 0.98 & \textbf{10\,(100\%)} & 0.98 \\
$200$ & 10 & 9 & \textbf{10\,(100\%)} & 597 & 9\,(90\%) & 0.96 & \textbf{10\,(100\%)} & 0.99 & \textbf{10\,(100\%)} & 0.99 & \textbf{10\,(100\%)} & 0.98 \\
$300$ & 10 & 9 & \textbf{10\,(100\%)} & 1395 & 9\,(90\%) & \textbf{1.03} & \textbf{10\,(100\%)} & 1.00 & \textbf{10\,(100\%)} & 0.99 & \textbf{10\,(100\%)} & 0.99 \\
$500$ & 10 & 0 & \textbf{0\,(0\%)} & -- & \textbf{0\,(0\%)} & -- & \textbf{0\,(0\%)} & -- & \textbf{0\,(0\%)} & -- & \textbf{0\,(0\%)} & -- \\
\bottomrule
\end{tabular}}
\end{table}

\begin{table}[!htbp]
\centering
\caption{Final optimality gaps for SLR instances not solved by any method within the time limit. Bold values indicate the best performance within each scenario count.}
\label{tab:gapdiff-slr-plain}
\small
\begin{tabular}{crccc}
\toprule
$|\mathcal{S}|$ & $|\mathcal{J}^0|$ & $\bar{g}_{\mathrm{B\&C}}$ (\%) & $\overline{\Delta g}_\mathrm{R}$ & $\overline{\Delta g}_\mathrm{S}$ \\
\midrule
$50$-$500$ & 0 & -- & -- & -- \\
$750$ & 4 & 13.6 & -0.1 & \textbf{+0.0} \\
$1000$ & 9 & 17.1 & \textbf{+5.3} & -0.1 \\
\bottomrule
\end{tabular}
\end{table}

\begin{table}[!htbp]
\centering
\caption{Final optimality gaps for \slrp\ instances not solved by any method within the time limit. Bold values indicate the best performance within each scenario count, and ties are all highlighted.}
\label{tab:gapdiff-slr-routing}
\small
\begin{tabular}{crccccc}
\toprule
$|\mathcal{S}|$ & $|\mathcal{J}^0|$ & $\bar{g}_{\mathrm{B\&C}}$ (\%) & $\overline{\Delta g}_\mathrm{R}$ & $\overline{\Delta g}_\mathrm{F}$ & $\overline{\Delta g}_\mathrm{S}$ & $\overline{\Delta g}_\mathrm{C}$ \\
\midrule
$50$-$300$ & 0 & -- & -- & -- & -- & -- \\
$500$ & 10 & \textbf{12.6} & -7.8 & \textbf{0.0} & -0.1 & -0.1 \\
\bottomrule
\end{tabular}
\end{table}

\clearpage
\bibliographystyle{abbrvnat}
\bibliography{references,extra}

\end{document}